\documentclass[11pt]{amsart}
\usepackage{amssymb, amsmath, amsthm, ulem}
\usepackage[utf8]{inputenc}
\usepackage{graphicx}
\usepackage[hidelinks]{hyperref}
\usepackage{xcolor}
\usepackage{epstopdf}
\usepackage{cancel}
\usepackage{tikz}
\usepackage{booktabs}
\usepackage{float}
\usetikzlibrary{positioning,arrows.meta,calc,angles,quotes,patterns}
 \usepackage{enumitem}

\newcommand{\dd}{\,\mathrm d}

\usepackage[a4paper, centering]{geometry}

\usepackage{scalerel,stackengine}

\usepackage{mathtools}
\usepackage{mathrsfs}

\newtheorem{theorem}{Theorem}
\newtheorem{proposition}[theorem]{Proposition}
\theoremstyle{plain}

\newtheorem{lemma}[theorem]{Lemma}
\newtheorem{corollary}[theorem]{Corollary}

\theoremstyle{definition}
\newtheorem{definition}[theorem]{Definition}

\theoremstyle{remark}
\newtheorem{remark}[theorem]{Remark}

\newcommand{\bt}{{\bf t}}

\newcommand{\one}{\mathbf 1}
\newcommand{\T}{\tau}
\providecommand{\Rcal}{\mathcal R}
\newcommand{\conv}{\operatorname{conv}}
\newcommand{\R}{\mathbb R}

\DeclareMathOperator{\Arg}{Arg}

\stackMath
\newcommand\reallywidecheck[1]{%
\savestack{\tmpbox}{\stretchto{%
  \scaleto{%
    \scalerel*[\widthof{\ensuremath{#1}}]{\kern-.6pt\bigwedge\kern-.6pt}%
    {\rule[-\textheight/2]{1ex}{\textheight}}
  }{\textheight}%
}{0.5ex}}%
\stackon[1pt]{#1}{\scalebox{-1}{\tmpbox}}%
}

\def\T{\tau}

\def\R{\mathbb{R}}

\newcommand{\Om}{\Omega}

\definecolor{proofblue}{HTML}{D6E9F3}
\definecolor{proofgreen}{HTML}{C4DFC8}
\definecolor{proofgold}{HTML}{F5DFAF}
\definecolor{proofink}{HTML}{172B3A}
\definecolor{figBlueDust}{HTML}{D7D7D7}
\definecolor{figBlueMid}{HTML}{DFE6EB}
\definecolor{figBluePale}{HTML}{F0F0F0}
\tikzset{
  psFill/.style={fill=figBlueDust},
  psHull/.style={fill=figBluePale},
  psOutline/.style={draw=black,line width=.75pt},
  psInternal/.style={draw=black!70,line width=.65pt},
  psDashed/.style={draw=black,dash pattern=on 2.4pt off 1.9pt,line width=.65pt}
}

\usepackage{float}

\providecommand{\conv}{\operatorname{conv}}

\providecommand{\Rcal}{\mathcal R}

\begin{document}
\title[Proof of the polygonal Saint-Venant inequality]{
Proof of the polygonal Saint-Venant inequality \\
conjectured by P\'olya--Szeg\"o 
}
\author[]{Beniamin Bogosel, Dorin Bucur, Ilaria Fragal\`a}

\thanks{}

\address[Beniamin Bogosel]{
	Faculty of Exact Sciences, Aurel Vlaicu University of Arad, 2 Elena Dr\u agoi Street, Arad, Romania}
\email {beniamin.bogosel@uav.ro}

\address[Dorin Bucur]{
Universit\'e  Savoie Mont Blanc, Laboratoire de Math\'ematiques CNRS UMR 5127 \\
  Campus Scientifique \\
73376 Le-Bourget-Du-Lac (France)
}
\email{dorin.bucur@univ-savoie.fr}

\address[Ilaria Fragal\`a]{
Dipartimento di Matematica \\ Politecnico  di Milano \\
Piazza Leonardo da Vinci, 32 \\
20133 Milano (Italy)
}
\email{ilaria.fragala@polimi.it}

\date{\today}


\begin{abstract}
We prove the polygonal Saint--Venant inequality conjectured by P\'olya and Szeg\"o: 
for every $n\geq 3$, among (possibly non--convex) polygons of prescribed area with at most $n$ sides, 
the regular $n$-gon uniquely maximizes torsional rigidity.
The proof is purely analytic and stems from the local version of the result given in our previous paper \cite{BBF}: it combines an area--angle proportional partition of a polygon, inspired by Solynin and Zalgaller \cite{SZ}, 
with a sharp comparison for the mixed torsional rigidity of the partition cells, and with 
an angular concavity property of this functional on the class of isosceles triangles. \end{abstract}
\maketitle

\section{Introduction and statement of the result}\label{sec:introduction}

The classical Saint--Venant inequality states that, among sets of prescribed area, torsional rigidity is maximized by the disc.
Recall that, 
 for any open bounded domain $\Omega$, 
  its torsional rigidity is given by 
\begin{equation}\label{eq:torsion-variational}
\T(\Omega ):=-\inf_{u\in H_0^1(\Omega)}\int_\Omega\bigl(|\nabla u|^2-2u\bigr)\,dx\,.
\end{equation}

In 1951, P\'olya and Szeg\"o conjectured that, for every $n\geq 3$, among polygons with $n$ sides and prescribed area, the regular $n$-gon has the largest torsional rigidity (see \cite[page 158]{PS}).
 They also proved the conjecture for $n=3,4$, the only cases in which Steiner symmetrization 
can be successfully applied as it preserves the number of sides. For  every $n \geq 5 $, the conjecture 
has remained open.   

Over the years, the main advances concerned related polygonal inequalities
or restricted families rather than the general conjecture. In \cite{FrGaLa}, 
the   extremality of the regular $n$--gon was proved in a subclass of convex polygons satisfying  an
asymmetry condition.  
More recently, Gui, Hu, Li and Zhang established strict monotonicity of
torsional rigidity along the sequence of equal-area regular polygons
\cite{GHLZ}, while Gui, Hu and Li proved a quantitative inequality for
tangential polygons \cite{GHL}.
A hybrid approach to the related polygonal Faber--Krahn problem was developed by the first and second authors, who showed that, for every fixed $n\geq 5$, the conjecture can be reduced to finitely many rigorous numerical verifications \cite{BB}. Building on this methodology, the strict local maximality of the regular $n$-gon for torsional rigidity, for $5\leq n\leq 25$,  was proved 
by the first author by means of a computer-assisted analysis based on validated estimates of the shape Hessian \cite{Bogosel26}.

In \cite{BBF}, we proved the local form of the conjecture for convex polygons: for every $n\geq5$, the regular $n$-gon is, up to rigid motions, the unique maximizer among convex $n$-gons sufficiently close to it.

In this work, we prove the conjecture in its full generality. 
Throughout the paper, a polygon is the bounded open region enclosed by a simple polygonal curve;  an $n$-gon is a polygon with exactly $n$ sides.

\begin{theorem}\label{thm:convex-optimality}
 Let $n\geq3$, let $D$ be a polygon with at most $n$ sides, and let $D_*$ be a regular $n$-gon of the same area. Then 
\[
\T(D)\leq\T(D_*),
\]
with equality if and only if $D$ is congruent to $D_*$.
\end{theorem}

The proof is entirely analytic. It starts with the construction of an area-angle proportional partition of $D$, inspired
 by the seminal paper \cite{SZ} by Solynin and Zalgaller
(who settled the analogous problem for logarithmic capacity), with a key difference: we construct a partition into convex polygons, rather than a covering by triangles.
The next step in the proof is a sharp comparison  result for the mixed torsional rigidity of the partition's cells. The proof concludes 
via an angular concavity result for the mixed torsional rigidity of isosceles triangles from our previous paper \cite{BBF}.  

\bigskip 
{\bf Plan of the paper.} The contents are organized as follows:

In Section~\ref{sec:convex} we explain the architecture of the proof of Theorem~\ref{thm:convex-optimality}: we state the required intermediate results and show how the theorem follows from them. 

The intermediate results are then proved in the subsequent sections:
Theorem \ref{thm:weighted} is proved in Section \ref{sec:partition}; 
Theorem \ref{lem:terminal} in Section \ref{sec:enhanced}; 
Theorem \ref{thm:convex-cap-comparison} in Section \ref{sec:convex-cap-proof}; 
Theorem \ref{thm:sectorial-comparison} in Section \ref{sec:spread}.

\bigskip 
{\bf Statement on AI assistance.} The proof  given   in this paper  is essentially the work of the authors, who
 take full responsibility for its whole content.  Throughout the manuscript, 
AI tools were used  to assist in checking 
mathematical arguments, suggesting some refinements, and clearing the exposition.

\bigskip
{\bf Addendum}:  Concerning the analogue problem  for the first Dirichlet eigenvalue, we point out that: 

\begin{itemize}

\item[(i)]  Adapting our approach, with the help of LLM models, 
gives a proof of the polygonal Faber--Krahn inequality. However, such a 
proof does not meet our objectives, because it contains hardly verifiable parts, 
it is lacking of  mathematical insight, and also relies on several supporting computational files.

\smallskip 
\item[(ii)]  The related preprint \cite{Wang26} appeared while submitting our manuscript. 
 
\end{itemize}

As of today, we believe that obtaining a transparent, human-verifiable proof of the polygonal Faber--Krahn inequality remains a worthwhile goal.

\section{Proof outline of Theorem~\ref{thm:convex-optimality} }\label{sec:convex}

In this section we present the proof strategy for Theorem~\ref{thm:convex-optimality}: we state its main ingredients and show how they combine to prove the theorem.

\subsection{Area--angle proportional partition.}  

\begin{theorem}[Area--angle proportional  partition of the convex hull]\label{thm:weighted}
Let $D$ be an $n$-gon, and let $A_1,\ldots,A_m$ be the
vertices of its convex hull $K=\operatorname{conv}(D)$,
ordered counterclockwise. 
There exist
$m$ pairwise disjoint open convex polygons $Q_1,\ldots,Q_m$, such that,  using indices modulo $m$:   
\begin{itemize}
\item{} for every $i = 1, \dots, m$, $[A_i, A_{i+1} ] \subseteq \partial Q_i$ 
\item{}  the sets $\overline Q_i$ form a partition of $\overline K$, i.e.
 \[
 \overline K= 
\bigcup_{i = 1, \dots, m}\overline{Q_i}\,, 
\]
\item{} for $i\in\mathcal I_+:=\{j:|Q_j\cap D|>0\}$, 
 $\overline Q_i\subseteq\overline T_i$ for a triangle $T_i$
with base $[A_i,A_{i+1}]$ and opposite angle 
\[
\alpha_i=2\pi \frac{|Q_i\cap D|}{|D|} \,. 
\]
\end{itemize} 
\end{theorem}
\begin{proof}
See Section~\ref{sec:partition}.
\end{proof}

\begin{remark} In the above theorem: 
\begin{itemize}
\item[(i)] When $D$ is convex, $K=D$, $m = n$ and $\mathcal I _+ = \{ 1, \dots, n \}$. 

\item[(ii)] The terminology {\it proportional } partition is  inspired from the seminal paper \cite{SZ}: it refers to the property that the quotient  between $|Q_i \cap D|$ and the angle $\alpha _i$ is independent of $i$. 
 
\item[(iii)] By construction,   the angles $\alpha_i$  have sum $2 \pi$ and the sets $\overline {Q_i \cap D}$ form a partition of $\overline D$: \begin{equation}\label{f:proppar} 
\sum_{i\in\mathcal I_+}\alpha_i=2\pi,\qquad
\overline D = \bigcup_{i\in\mathcal I_+}\overline{Q_i \cap D}.
\end{equation} 
\end{itemize} 
\end{remark} 

\begingroup
\begin{figure}[H]
\centering
\begin{minipage}[t]{0.485\linewidth}
\centering
\begin{tikzpicture}[x=1.76cm,y=1.76cm,line join=round,line cap=round,font=\scriptsize,every node/.append style={text=black}]
  \coordinate (A1) at (1.092512,0.663366);
  \coordinate (A2) at (-1.092512,0.663366);
  \coordinate (A3) at (-1.769852,-0.578146);
  \coordinate (A4) at (0.000000,-1.081900);
  \coordinate (A5) at (1.769852,-0.578146);
  \fill[psHull] (-0.296087,-0.130193)--(-1.092512,0.663366)--(1.092512,0.663366)--(0.296087,-0.130193)--(-0.000000,-0.170457)--cycle;
  \fill[psHull] (-1.769852,-0.578146)--(-1.092512,0.663366)--(-0.296087,-0.130193)--cycle;
  \fill[psHull] (0.000000,-1.081900)--(-1.769852,-0.578146)--(-0.296087,-0.130193)--(0.000000,-0.170457)--cycle;
  \fill[psHull] (0.000000,-1.081900)--(0.000000,-0.170457)--(0.296087,-0.130193)--(1.769852,-0.578146)--cycle;
  \fill[psHull] (1.769852,-0.578146)--(0.296087,-0.130193)--(1.092512,0.663366)--cycle;
  \draw[psInternal] (-0.296087,-0.130193)--(-1.092512,0.663366);
  \draw[psInternal] (-0.000000,-0.170457)--(-0.296087,-0.130193);
  \draw[psInternal] (0.296087,-0.130193)--(-0.000000,-0.170457);
  \draw[psInternal] (1.092512,0.663366)--(0.296087,-0.130193);
  \draw[psInternal] (-0.296087,-0.130193)--(-1.769852,-0.578146);
  \draw[psInternal] (0.000000,-0.170457)--(0.000000,-1.081900);
  \draw[psInternal] (0.296087,-0.130193)--(1.769852,-0.578146);
  \draw[psDashed] (A1)--(-0.000000,-0.425215);
  \draw[psDashed] (A2)--(-0.000000,-0.425215);
  \draw[psDashed] (A3)--(0.000000,-0.040197);
  \draw[psDashed] (A4)--(0.000000,-0.040197);
  \draw[psDashed] (A5)--(0.000000,-0.040197);
  \draw[psOutline] (A1)--(A2)--(A3)--(A4)--(A5)--cycle;
  \node[inner sep=.30pt] at (0.000000,0.390000) {$Q_1$};
  \node[inner sep=.30pt] at (-1.265000,-0.100000) {$Q_2$};
  \node[inner sep=.30pt] at (-0.800000,-0.600000) {$Q_3$};
  \node[inner sep=.30pt] at (0.800000,-0.600000) {$Q_4$};
  \node[inner sep=.30pt] at (1.265000,-0.100000) {$Q_5$};
  \node[above right,inner sep=1.0pt] at (A1) {$A_1$};
  \node[above left,inner sep=1.0pt] at (A2) {$A_2$};
  \node[below left,inner sep=1.0pt] at (A3) {$A_3$};
  \node[below,inner sep=1.0pt] at (A4) {$A_4$};
  \node[below right,inner sep=1.0pt] at (A5) {$A_5$};
\end{tikzpicture}
\end{minipage}
\hfill
\begin{minipage}[t]{0.485\linewidth}
\centering
\begin{tikzpicture}[x=1.65cm,y=1.76cm,font=\scriptsize,line cap=round,line join=round,every node/.append style={text=black}]
  \coordinate (A1) at (1.092512,0.663366);
  \coordinate (A2) at (-1.092512,0.663366);
  \coordinate (A3) at (-1.769852,-0.578146);
  \coordinate (A4) at (0.000000,-1.081900);
  \coordinate (A5) at (1.769852,-0.578146);
  \coordinate (X2) at (-1.092512,-0.413462);
  \coordinate (X5) at (1.092512,-0.413462);
  \coordinate (Y) at (0.000000,-0.646922);
  \coordinate (H3) at (-1.485373,-0.508980);
  \coordinate (H4) at (1.485373,-0.508980);
  \coordinate (W) at (-1.329581,0.228837);
  \coordinate (Z2) at (-1.300000,-0.120000);
  \coordinate (Z5) at (1.348785,-0.068350);
  \coordinate (V) at (0.000000,-0.987182);
  \fill[psFill] (A1)--(A2)--(A3)--(A4)--(A5)--cycle;
  \fill[psHull] (V)--(H4)--(X5)--(A1)--(Z5)--(A5)--(A4)--(A3)--(Z2)--(W)--(A2)--(X2)--(H3)--cycle;
  \draw[psInternal] (A2)--(X2)--(Y)--(X5)--(A1);
  \draw[psInternal] (A3)--(X2);
  \draw[psInternal] (A4)--(Y);
  \draw[psInternal] (A5)--(X5);
  \draw[psInternal] (A1)--(A2)--(A3)--(A4)--(A5)--cycle;
  \draw[psDashed,
    preaction={draw=white,line width=1.35pt}] (A2)--(A3)--(X2)--cycle;
  \draw[psDashed,
    preaction={draw=white,line width=1.35pt}] (A3)--(A4)--(0.000000,-0.147835)--cycle;
  \draw[psDashed,
    preaction={draw=white,line width=1.35pt}] (A4)--(A5)--(0.000000,-0.147835)--cycle;
  \draw[psDashed,
    preaction={draw=white,line width=1.35pt}] (A5)--(A1)--(X5)--cycle;
  \draw[psOutline] (V)--(H4)--(X5)--(A1)--(Z5)--(A5)--(A4)--(A3)--(Z2)--(W)--(A2)--(X2)--(H3)--cycle;
  \node[inner sep=0pt] at (0.00,0.34) {$Q_1$};
  \node[inner sep=0pt] at (-1.30,-0.33) {$Q_2$};
  \node[inner sep=0pt] at (-0.70,-0.64) {$Q_3$};
  \node[inner sep=0pt] at (0.70,-0.64) {$Q_4$};
  \node[inner sep=0pt] at (1.33,-0.33) {$Q_5$};
  \node[inner sep=0pt] at (-1.87,-0.33) {$T_2$};
  \node[inner sep=0pt] at (-0.75,-1.02) {$T_3$};
  \node[inner sep=0pt] at (0.75,-1.02) {$T_4$};
  \node[inner sep=0pt] at (1.87,-0.33) {$T_5$};
  \node[above right,inner sep=1pt] at (A1) {$A_1$};
  \node[above left,inner sep=1pt] at (A2) {$A_2$};
  \node[below left,inner sep=1pt] at (A3) {$A_3$};
  \node[below,inner sep=1pt] at (A4) {$A_4$};
  \node[below right,inner sep=1pt] at (A5) {$A_5$};
\end{tikzpicture}
\end{minipage}
\caption{Two examples of proportional convex hull partitions.}
\label{fig:weighted-convex-and-nonconvex}
\end{figure}
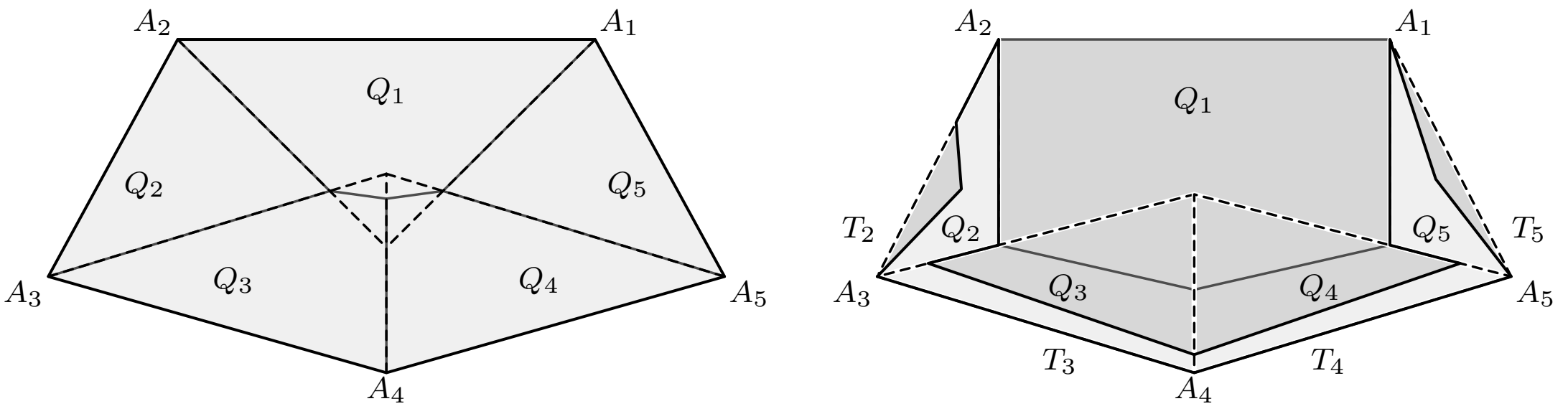

\endgroup

When $D$ is nonconvex, some vertices of $D$ may lie inside the
cells $Q_i$ or in the relative interior of their bases $[A_i,A_{i+1}]$ (as occurs for $Q_2$ and $Q_5$ in Figure~\ref{fig:weighted-convex-and-nonconvex}). The next
theorem shows that, starting from the partition of Theorem~\ref{thm:weighted},   it is possible to  obtain a refinement of the partition of $D$ in which
both situations are avoided and the crucial properties~\eqref{f:proppar} are preserved.
This is achieved by replacing the pairs $\{(Q_i,T_i)\}_{i\in\mathcal I_+}$ with
at most $n$ sectorial pairs $\{(Q_j,\Delta_j)\}$.  
To be more precise, let us introduce the following:

\begin{definition}\label{def:sectorial}
A {\it sectorial pair} is a pair $(Q,\Delta)$, where $\Delta$ is a triangle
with opening $\theta\in(0,\pi)$ and opposite base $\Gamma_\theta$,
and $Q\subseteq\Delta$ is a convex polygon with
$\Gamma_\theta\subseteq\partial Q$.
The pair is {\it simple with respect to $D$} if no vertex of $D$ belongs to
$Q\cup\operatorname{relint}\Gamma_\theta$.
The angle and base are called the {\it opening} and {\it base} of the pair. 
\end{definition}

\begin{theorem}[Enhanced area-angle proportional partition of an $n$-gon]\label{lem:terminal}
Let $D$ be an $n$-gon and let $K:=\operatorname{conv}(D)$.   There exist at most $n$ 
sectorial pairs $(Q _j, \Delta _j)$, with openings $\theta _j$ and bases $\Gamma _{\theta _j} \subset  \partial K$, which are
simple with respect to $D$, and  such that the  (possibly disconnected) sets
 $U _j := D \cap Q _j$ enjoy the following properties: 
 
 \smallskip 
\begin{itemize}
\item{} $U _j$ have positive measures, are pairwise disjoint 
and their closures form a partition of  $\overline D$, i.e.  
$$\overline D = \bigcup_{j } \overline U _j\,.$$    
\smallskip 

\item{} the above partition is area--angle proportional, in the sense that
$$\theta _j = 2 \pi \frac{|U _j |} {|D|  }\,.$$  
\end{itemize}

\end{theorem}

\begin{proof}
See Section~\ref{sec:enhanced}.
\end{proof}
 
\begin{remark} By construction, the angles $\theta_j$ in the above theorem satisfy
 \begin{equation}\label{f:proppar2}
\sum_{j} \theta_j
=2\pi.
\end{equation}
\end{remark}

\subsection{Optimization of the mixed torsional rigidity} 

\smallskip \leavevmode Let $\Om$ be a bounded open set with finitely many Lipschitz connected components $\Om_h$, and let $\Gamma\subset\partial\Om$ be such that each $\Gamma_h:=\Gamma\cap\partial\Om_h$ has positive length. All trace and normal conditions below are understood componentwise. The complementary portion $\partial \Om \setminus\Gamma$ is the Neumann boundary. We set
\begin{equation}\label{eq:space-V}
H_\Gamma(\Om):=\bigl\{u\in H^1(\Om):u=0\text{ on }\Gamma\bigr\}.
\end{equation}

\noindent The {\it mixed torsional rigidity} of $\Om$ with Dirichlet portion $\Gamma$ is defined, in its primal formulation, by
\begin{equation}\label{eq:mixed-primal}
\T_{\rm mix}(\Om,\Gamma)=-\inf_{u\in H_\Gamma(\Om)}
\int_{ \Om} \bigl(|\nabla u|^2-2u\bigr)\,dx=\int_\Om  w_{\Om,\Gamma}\,dx,
\end{equation}
\par\noindent where $w_{\Om,\Gamma}$ solves
\[
\begin{cases}
-\Delta w=1 & \text{in $\Om$},\\
w=0 & \text{on }\Gamma,\\
\partial_\nu w=0 & \text{on }\partial \Om\setminus\Gamma.
\end{cases}
\]
\par\noindent The function $w_{\Om,\Gamma}$ is called the {\it mixed torsion function of $\Om$ (with Dirichlet portion $\Gamma$)}.

\noindent It will be useful to recall that $\T_{\rm mix}(\Om,\Gamma)$ can be recast in   dual formulation as
\[
\T_{\rm mix}(\Om,\Gamma)=\min_{\sigma\in X_\Gamma(\Om)}\int_ \Om |\sigma|^2\,dx.
\]

\begingroup
\noindent where 
\begin{equation}\label{eq:space-H}
X_\Gamma(\Om):=\bigl\{\sigma\in L^2(\Om ; \mathbb R^2):
\operatorname{div}\sigma=-1\text{ in } \Om ,\quad
\sigma\cdot\nu=0\text{ on }\partial  \Om \setminus\Gamma\bigr\}.
\end{equation}
Here the divergence is understood in the distributional sense, and the normal trace is understood in the weak sense.
The elements of $X_\Gamma(\Om)$ are called {\it feasible fluxes}.
\endgroup

\noindent The minimum in the dual formulation is attained at $\sigma=\nabla w_{\Om,\Gamma}$.

We observe that, if $\Om$ is disconnected  with connected components $\Om _h$, letting $\Gamma _h:= \Gamma \cap \partial \Om _h$, there holds  
\[
\T_{\rm mix}(\Om,\Gamma)=\sum_h\T_{\rm mix}(\Om_h,\Gamma_h).
\]
 When the Dirichlet portion $\Gamma$ will be clear from the context, it will be dropped and we shall use for simplicity the notation 
 $\T_{\rm mix} (\Om)$  and $w _\Om$.

It is convenient to introduce also the following area-rescaled mixed torsional rigidity functional, 
which is invariant under domain dilation:
  
  \begin{equation}\label{f:rescaled} \Rcal(\Om,\Gamma)=\frac{\T_{\rm mix}(\Om,\Gamma) }{|\Om|^2}  \,,
  \end{equation} 
For any given $\theta \in (0, \pi)$, we  denote by $I _\theta$ an isosceles triangle with apex angle $\theta$. 

\smallskip
{\it For brevity, we shall denote by $\mathcal R (\theta)$
the rescaled mixed torsional rigidity functional  defined by  \eqref{f:rescaled} when $\Om = I _\theta$ and $\Gamma$ is the side of $I _\theta$ opposite to $\theta$. } 

\smallskip 

We are going to state two optimization results  for the rescaled mixed torsional rigidity.

In the first optimization result,  we fix a triangle  $T$  with an interior angle $\theta\in(0,\pi)$ and opposite side $\Gamma _\theta$,  and 
 we maximize $\mathcal R$ over the following
class of convex {polygonal} caps: 
 \begin{equation}\label{f:capspol} 
\mathcal C_\theta(T):=\bigl\{Q: Q\text{ is a convex polygon},\ |Q|>0,\ Q\subseteq T,\ \Gamma_\theta\subseteq\partial Q\bigr\}\,.
\end{equation}

\begin{theorem}[Optimization for convex {polygonal} caps in a triangle]\label{thm:convex-cap-comparison}
\smallskip 
Let $T$ be a given triangle with a distinguished interior angle $\theta\in(0,\pi)$, and opposite side $\Gamma _\theta$.  Let $\mathcal C_\theta(T)$ be the class 
defined in \eqref{f:capspol}. For every $Q\in\mathcal C_\theta(T)$ it holds
  \begin{equation}\label{eq:convex-cap-comparison}
\Rcal(Q,\Gamma_\theta)\le \Rcal(\theta).
\end{equation}
Equality holds if and only if $Q=T$ and $T$ is similar to $I_\theta$.
\end{theorem}

\begin{proof}
See Section~\ref{sec:convex-cap-proof}.
\end{proof}

\medskip
 The second optimization result, whose proof relies on Theorem \ref{thm:convex-cap-comparison},  concerns the rescaled mixed torsional rigidity of the sets $U _j$ issued from the  enhanced proportional partition of a $n$-gon provided by Theorem \ref{lem:terminal}. 
\begin{theorem}[Optimization for polygonal sets in a sectorial pair] \label{thm:sectorial-comparison} Let $D$ be a $n$-gon, and let 
$(Q _j, \Delta _j)$ be sectorial pairs of openings $\theta _j$, 
simple with respect to $D$,  as given by Theorem \ref{lem:terminal}, 
and let  
$U _j := D \cap Q _j$. Then for every $j$ it holds  
$$\mathcal R ( U _ j,  \partial U _ j \cap \partial D) \leq \mathcal R (\theta_j) \,.$$ 
Moreover, equality holds if and only if
$U_j$ is an isosceles triangle with apex angle $\theta_j$
and base $\partial U_j\cap\partial D$.

\end{theorem} 

\proof See Section \ref{sec:spread}. \qed

\subsection{Conclusion of the proof of Theorem~\ref{thm:convex-optimality} } 

 Before proceeding, we recall the following result for the function  
\begin{equation}\label{eq:triangular-F}
 \mathcal F(\theta):=\theta^2\Rcal(I_\theta),\qquad \theta\in(0,\pi).
\end{equation}

\begin{lemma}\label{cc:lem:triangle-benchmark}
The function $\mathcal F$ is smooth and satisfies
\begin{equation}\label{f:effe} 
 \mathcal F''(\theta)<0\qquad\text{for every }\theta\in(0,\pi).
\end{equation} 
\end{lemma}

\begin{proof}
This follows from~\cite[Theorem~8]{BBF}, where the function
$F=\mathcal F/4$ is considered.
\end{proof}

\bigskip

{\it  Proof of Theorem~\ref{thm:convex-optimality}.}  
Assume that $n \geq 5$ (otherwise Theorem~\ref{thm:convex-optimality} is already known), and  assume also, without loss of generality, that $D$ is a $n$-gon (i.e., it has exactly $n$ sides). 
Let $S=|D|$. By Theorem~\ref{lem:terminal} there exist
$n$  sectorial pairs $(Q_j,\Delta_j)$, with openings
$\theta_j$, such that, setting $U_j=D\cap Q_j$,
\[
|U_j|=\frac{S\theta_j}{2\pi},
\qquad
\sum_{j=1}^k\theta_j=2\pi.
\]
Set $\Gamma_j=\partial U_j\cap\partial D$.
We have
\[
\begin{aligned}
\frac{4\pi^2}{S^2}\T(D)
&\underset{\text{\scriptsize(1)}}{\le}
\frac{4\pi^2}{S^2}\sum_{j=1}^k\T_{\rm mix}(U_j,\Gamma_j)
\underset{\text{\scriptsize(2)}}{=}
\sum_{j=1}^k\theta_j^2\Rcal(U_j,\Gamma_j)\\
&\underset{\text{\scriptsize(3)}}{\le}
\sum_{j=1}^k\theta_j^2\Rcal(\theta_j)
\underset{\text{\scriptsize(4)}}{=}
\sum_{j=1}^k\mathcal F(\theta_j)
\underset{\text{\scriptsize(5)}}{\le}
k\mathcal F\left(\frac{2\pi}{k}\right)\\
&\underset{\text{\scriptsize(6)}}{\le}
n\mathcal F\left(\frac{2\pi}{n}\right)
\underset{\text{\scriptsize(7)}}{=}
\frac{4\pi^2}{S^2}\T(D_*).
\end{aligned}
\]
 Here (1) follows from Neumann relaxation across the
interfaces of the partition; (2) from the definition
of $\Rcal$ and the proportionality of the partition;
(3) from Theorem~\ref{thm:sectorial-comparison};
(4) from the definition of $\mathcal F$;
(5) from Jensen's inequality and $\sum_j\theta_j=2\pi$;
(6) from the strict increasing monotonicity of
$t\mapsto t\mathcal F(2\pi/t)$ for $t>2$, which holds as a consequence
of strict concavity of $\mathcal F$ and of the equality $\mathcal F(0+)=0$ (see the proof of Corollary~\ref{c:monobench});
(7) from the decomposition of $D_*$ into $n$ congruent
isosceles triangles along its symmetry axes.

If equality holds, strict concavity and strict
monotonicity imply
\[
k=n \qquad  \text{ and } \qquad  \theta_j=\frac{2\pi}{n}
\quad \forall j=1,\ldots,n,
\]
while equality in (3) gives
\[
\Rcal(U_j,\Gamma_j)=\Rcal(\theta_j)
\qquad  \forall j=1,\ldots,n.
\]
By the equality case of Theorem~\ref{thm:sectorial-comparison},
each $U_j$ is an isosceles triangle with apex angle $2\pi/n$.
Since all have area $S/n$, they are congruent.
Each triangular base lies on $\partial D$; equality in the
Neumann relaxation excludes any other boundary side of positive
length on $\partial D$. Hence the $n$ bases cover all $n$  
sides of $D$. A straight base lies in a single side, so counting
shows that each base is one full side. At every vertex of $D$,
the two incident isosceles triangles contribute base angles
$\frac12(\pi-2\pi/n)$ each; thus the interior angle of $D$
is at least $\pi-2\pi/n$. Their sum is $(n-2)\pi$, so all
these inequalities are equalities. The polygon is therefore
convex and equiangular; its full sides, being bases of congruent
triangles, have equal length. Hence $D$ is regular. 
\qed

\bigskip

\section{Proof of Theorem~\ref{thm:weighted}}\label{sec:partition}

We denote by $h_i$ the signed Euclidean distance to the supporting line of $[A_i, A_{i+1}]$, nonnegative on $\overline K$ and vanishing there exactly on $[A_i, A_{i+1}]$.

For $\bt=(t_1,\dots,t_m)\in\R^m$  put
\begin{equation}\label{eq:cells}
 Q_i(\bt)=\left\{x\in K:e^{t_i}h_i(x)<e^{t_j}h_j(x)
 \text{ for all }j\ne i\right\}.
\end{equation}
The open convex cells $Q_i(\bt)$ are pairwise disjoint, and their
closures cover $\overline K$. Since $h_i=0$ on $[A_i, A_{i+1}]$ and the other
$h_j$ are positive on its relative interior,
$[A_i, A_{i+1}]\subseteq\partial Q_i(\bt)$ and $Q_i(\bt)\ne\varnothing$.
Replacing $\bt$ by  $\bt+c\one$ does not change the geometry of the partition.

At $A_i$ the dividing ray of $Q_{i-1}$ and $Q_i$ makes an angle
$x_i(\bt)\in(0,\beta_i)$ with the outgoing side $[A_i, A_{i+1}]$. Its equation is
\[
 e^{t_i-t_{i-1}}\sin x_i=\sin(\beta_i-x_i)\,,
\]

where 
\[
\beta_i:=\angle A_{i-1}A_iA_{i+1}\in(0,\pi).
\]

Define
\begin{equation}\label{eq:phi}
 \varphi_i(s)=\Arg\left(e^s+\cos\beta_i+
  \mathrm i\sin\beta_i\right),\qquad
 \varphi_i'(s)=-\frac{e^s\sin\beta_i}
 {e^{2s}+2e^s\cos\beta_i+1}<0.
\end{equation}
Then $x_i(\bt)=\varphi_i(t_i-t_{i-1})$,
$\varphi_i(-\infty)=\beta_i$ and $\varphi_i(+\infty)=0$.
The \textit{formal} angle opposite $[A_i, A_{i+1}]$ is
\begin{equation*} 
 a_i(\bt)=\pi-\beta_{i+1}-x_i(\bt)+x_{i+1}(\bt).
\end{equation*}
For general weights it can be negative; nevertheless
\begin{equation*} 
 \sum_{i=1}^m a_i(\bt)=m\pi-\sum_i\beta_i=2\pi,
 \qquad a_i(\bt)<\pi.
\end{equation*}
We write
  $$
 \Phi_i(s)=\int_0^s\varphi_i(u)\dd u,$$
 $$\mathscr A(\bt)=\sum_i(\pi-\beta_{i+1})t_i
 -\sum_i\Phi_i(t_i-t_{i-1}), $$
 \begin{equation} \label{eq:areapotential}
 \mathscr B(\bt)=\int_D \min_i\left(t_i+\log\frac{h_i(x)}{L}\right)\dd x,
 \qquad L>0.
\end{equation}
The logarithmic singularities on the supporting lines are integrable.
The minimum in \eqref{eq:areapotential} is an almost-everywhere uniquely
attained concave function of $\bt$, since ties between distinct positive
weighted affine distances occur on lines. Consequently,
\begin{equation}\label{eq:gradients}
 \frac{\partial\mathscr A}{\partial t_i}=a_i(\bt),
 \qquad \frac{\partial\mathscr B}{\partial t_i}=|Q_i(\bt)\cap D|.
\end{equation}
In particular $\mathscr B$ is $C^1$ and concave. Differentiation and
\eqref{eq:phi} give
\begin{equation*} 
 v^T\nabla^2\mathscr A(\bt)v
 =\sum_i[-\varphi'_i(t_i-t_{i-1})](v_i-v_{i-1})^2.
\end{equation*}
Every coefficient is positive. Thus $\mathscr A$ is strictly convex
modulo the constant vectors. Each $\Phi_i$ is bounded above because
$\Phi_i(s)\le0$ for $s\le0$ and $\varphi_i(s)=O(e^{-s})$ for $s\to+\infty$.

Set
\begin{equation*}
 \lambda=\frac{2\pi}{|D|},\qquad
 \Psi(\bt)=\mathscr A(\bt)-\lambda\mathscr B(\bt).
\end{equation*}
It is convex, differentiable and strictly convex modulo constants.
The shift laws
$\mathscr A(\bt+c\one)=\mathscr A(\bt)+2\pi c$ and
$\mathscr B(\bt+c\one)=\mathscr B(\bt)+c|D|$
make $\Psi$ invariant under common translations.
Normalize by $\min_i t_i=0$. All $t_i\ge0$; because the
$\Phi_i$ are uniformly bounded above and
$c_0:=\min_i(\pi-\beta_{i+1})>0$,
\begin{equation}\label{eq:coercive1}
 \mathscr A(\bt)\ge c_0\sum_i t_i-C_0.
\end{equation}
If $t_j=0$, then
\begin{equation*} 
 \mathscr B(\bt)\le\int_D\log\frac{h_j(x)}L\dd x\le M,
\end{equation*}
with $M$ independent of $\bt$ (only finitely many $j$).
Thus $\Psi$ is coercive on this closed gauge set, attains a minimum
$\bt^*$, and its invariance makes this a global unconstrained minimum.
By \eqref{eq:gradients},
\begin{equation*} 
 0=\partial_{t_i}\Psi(\bt^*)
 =a_i(\bt^*)-\frac{2\pi}{|D|}|Q_i(\bt^*)\cap D|.
\end{equation*}
Strict convexity on the quotient yields uniqueness modulo constants.

Set $\alpha_i=a_i(\bt^*)$ and $Q_i=Q_i(\bt^*)$.
If $\alpha_i>0$, the adjacent dividing rays issuing from $A_i$ and
$A_{i+1}$ have base angles $x_i>0$ and
$\beta_{i+1}-x_{i+1}>0$ summing to $\pi-\alpha_i<\pi$.
They meet forward at $V_i$, giving a triangle $T_i$ of opposite angle
$\alpha_i$. The two adjacent inequalities in \eqref{eq:cells}, together
with $h_i\ge0$, are the three inward halfplane conditions of $T_i$;
therefore $\overline Q_i\subseteq\overline T_i$.

Let $\mathcal I=\{i:\alpha_i>0\}$. If a point of $D$ were outside
the closed union $\bigcup_{i\in\mathcal I}\overline Q_i$,
a small open disk of $D$ would be covered, up to the separating line segments, by finitely many inactive cells,
all having zero area in $D$. This is impossible. Therefore
$\overline D\subseteq\bigcup_{i\in\mathcal I}\overline Q_i
\subseteq\bigcup_{i\in\mathcal I}\overline T_i$.
Theorem~\ref{thm:weighted} follows.
\qed

 \bigskip
 
\section{Proof of Theorem \ref{lem:terminal}}\label{sec:enhanced}

Let us consider the proportional partition of $D$ given by Theorem~\ref{thm:weighted}. We retain only the pairs $(Q_i,T_i)$ with $|Q_i\cap D|>0$.

Call a vertex of $D$ \textit{disruptive} if it lies in
$Q_i$ or in the relative interior
of the base of some retained pair $(Q_i,T_i)$. 

There are at most $n-m$ disruptive vertices, where $m$ is
the number of vertices of $K$. 

Let us fix a disruptive vertex $Z$, and let $Q\subset T:=\Delta(O,A,B)$
with opening $\theta$ and base $[A,B]\subset\partial K$ be such that
$Z\in Q\cup(A,B)$. 

We apply the following refinement procedure to make $Z$ no longer disruptive. 
For $t\in[0,\theta]$, draw through $Z$ the line parallel to the
direction making angle $t$ with $OA$ toward $OB$; it meets $[A,B]$
at $C_t$. The line divides $Q$ into convex portions $Q_\ell(t)$ and
$Q_r(t)$, with respective bases $[A,C_t]$ and $[C_t,B]$.
Each portion is contained in a triangle with the same base and
opposite angle $t$ or $\theta-t$, respectively: its lateral lines
are the cutting line and one lateral line of $T$.

The function
$f(t):=\theta|D\cap Q_\ell(t)|-t|D\cap Q|$ is continuous, with
$f(0)\ge0$ and $f(\theta)\le0$. Thus $f(t_*)=0$ for some
$t_*\in[0,\theta]$. Consequently, the two newly created cells satisfy
\[ 
|D\cap Q_\ell(t_*)|=\frac{t_*}{\theta}|D\cap Q|,
\qquad
|D\cap Q_r(t_*)|=\frac{\theta-t_*}{\theta}|D\cap Q|.
\]
If $0<t_*<\theta$, retain both positive-area cells, with
openings $t_*$ and $\theta-t_*$. If $t_*=0$ or $t_*=\theta$,
discard the zero-mass cell  and retain the other, with its
nondegenerate containing triangle of opening $\theta$.
In every case the retained cells have positive opening and
preserve area--angle proportionality. The vertex $Z$ lies on
the cut (or at an endpoint of the new bases), hence outside every
retained cell and the relative interior of every retained base.
No other vertex becomes disruptive. Therefore each refinement
reduces their total number by at least one and increases the
number of positive-area cells by at most one.

After at most $n-m$ refinements, all cells are simple, and their
number is at most $m+(n-m)=n$.
Moreover, the cells remain convex, with bases on $\partial K$,
their intersections $U_j$ with $D$ cover $D$ (up to a null set),
and the proportionality property
$\theta_j=2\pi|U_j|/|D|$ is preserved.

Therefore, the resulting family of at most $n$ simple cells,
together with their containing triangles, gives the asserted
family of sectorial pairs. \qed 
\bigskip

\section{\texorpdfstring{Proof of Theorem~\ref{thm:convex-cap-comparison}}{Proof of the convex-cap comparison}}\label{sec:convex-cap-proof}

{
Throughout this section, $T=\triangle(O,A,B)$ has apex angle
$\theta\in(0,\pi)$, and $Q\in\mathcal C_\theta(T)$ has Dirichlet
boundary $[A,B]$. We denote by $H_\Gamma(Q)$ and $X_\Gamma(Q)$ the spaces defined
in~\eqref{eq:space-V} and~\eqref{eq:space-H}, respectively, with
$\Gamma=[A,B]$. We write for brevity
\[
 \Rcal(\theta):=\Rcal(I_\theta).
\]
}

\subsection{\texorpdfstring{Variational bounds and the large-side-ratio regime}{Variational bounds and the large-side-ratio regime}}\label{cc:sec:bounds}

We are going to obtain some upper bounds for $\Rcal (Q)$ and some lower bounds for $\Rcal(\theta)$, by exploiting the 
variational characterization of the mixed torsional rigidity of a domain $D$ recalled in Section~\ref{sec:convex}. 
We use the dual formulation to obtain upper bounds for $\Rcal(Q)$ and the primal formulation to obtain lower bounds for $\Rcal(\theta)$. Thus, 
for every nonzero
$v\in H_\Gamma(D)$ and $\sigma\in X_\Gamma(D)$,
\begin{equation}\label{cc:eq:duality}
 \frac{(\int_D v)^2}{\int_D|\nabla v|^2}
 \le\T_{\rm mix}(D)\le\int_D|\sigma|^2.
\end{equation}

\begin{proposition}[Upper bounds for $\Rcal(Q)$]\label{cc:prop:bounds}
Set
\[
 t:=\max\left\{\frac{|OA|}{|OB|},\frac{|OB|}{|OA|}\right\},
 \qquad z:=\cos\theta,\qquad \delta:=\frac{|Q|}{|T|}.
\]
Then
\begin{align}
 \Rcal(Q)&\le\frac{\delta}{3t\sin\theta},\label{cc:eq:column}\\[\smallskipamount]
 {\Rcal(Q)}&{\le\frac{\delta\sin\theta}{3(t+t^{-1}-2z)}
 \qquad\text{if $O$ projects orthogonally inside $[A,B]$,}}\label{cc:eq:normal-roof}\\[\smallskipamount]
 {\Rcal(Q)}&{\le\frac\delta6\cot\frac\theta2,}\label{cc:eq:uniform-angle}
\end{align}
\end{proposition}
\begin{proof}
Let $a=|OA|>0$ and $b=|OB|>0$, and let $e_1,e_2$ be the unit vectors
pointing along the sides $[O,A]$ and $[O,B]$, respectively. Thus
$A=ae_1$ and $B=be_2$. We represent $Q$ as the image, under the
ray-coordinate map
$X(x,y):=ax e_1+by e_2$, of a convex  polygonal  cap $\widehat Q$ of the unit simplex
\[
 \Delta=\{x,y>0:x+y<1\}.
\]
More precisely,
\[
 \widehat Q=\{(x,y):0<x<1,\ \ell(x)<y<1-x\},
\]
where $\ell$ is nonnegative and convex. Hence
$h(x):=1-x-\ell(x)$ is nonnegative and concave. Since
$|\det DX|=ab\sin\theta$ and $|T|=ab\sin\theta/2$, we have
\[
 \int_0^1h(x)\,\dd x=\frac{|Q|}{ab\sin\theta}=\frac\delta2.
\]
Figure~\ref{cc:fig:clipped-cap} illustrates a convex {polygonal} cap and its preimage under $X$.
\par
\begin{figure}[htbp]
\centering
\begin{minipage}[c]{.47\linewidth}
\raggedleft
\begin{tikzpicture}[scale=.80,font=\small,line join=round,line cap=round,draw=black,every node/.append style={text=black}]
\coordinate (O) at (0,0);
\coordinate (A) at (5.5,0);
\coordinate (B) at (2.1,2.8);
\coordinate (P) at (1.1,0);
\coordinate (R) at (1.47,1.96);
\coordinate (U) at (.986,.728);
\fill[psHull] (P)--(A)--(B)--(R)--(U)--cycle;
\draw[psDashed] (O)--(P) (O)--(R);
\draw[psOutline] (P)--(A)--(B)--(R)--(U)--cycle;
\draw[psInternal] (.65,0) arc[start angle=0,end angle=53.13,radius=.65];
\node[below left=2pt] at (O) {$O$};
\node[below=3pt] at (A) {$A$};
\node[above=3pt] at (B) {$B$};
\node at (.82,.31) {$\theta$};
\node at (2.85,1.05) {$Q$};
\end{tikzpicture}
\end{minipage}\hspace{.04\linewidth}%
\begin{minipage}[c]{.45\linewidth}
\raggedright
\begin{tikzpicture}[scale=.72,font=\small,line join=round,line cap=round,draw=black,every node/.append style={text=black}]
\coordinate (O) at (0,0);
\coordinate (A) at (4.1,0);
\coordinate (B) at (0,4.1);
\coordinate (P) at (.82,0);
\coordinate (R) at (0,2.87);
\coordinate (U) at (.328,1.066);
\fill[psHull] (P)--(A)--(B)--(R)--(U)--cycle;
\draw[psDashed] (O)--(P) (O)--(R);
\draw[psOutline] (P)--(A)--(B)--(R)--(U)--cycle;
\node[below left=2pt] at (O) {$(0,0)$};
\node[below=3pt] at (A) {$(1,0)$};
\node[above=3pt] at (B) {$(0,1)$};
\node at (1.85,1.35) {$\widehat Q$};
\end{tikzpicture}
\end{minipage}
\caption{A convex {polygonal} cap $Q\in\mathcal C_\theta(T)$ (left) and its preimage
$\widehat Q$ under the ray-coordinate map $X$ (right).}
\label{cc:fig:clipped-cap}
\end{figure}
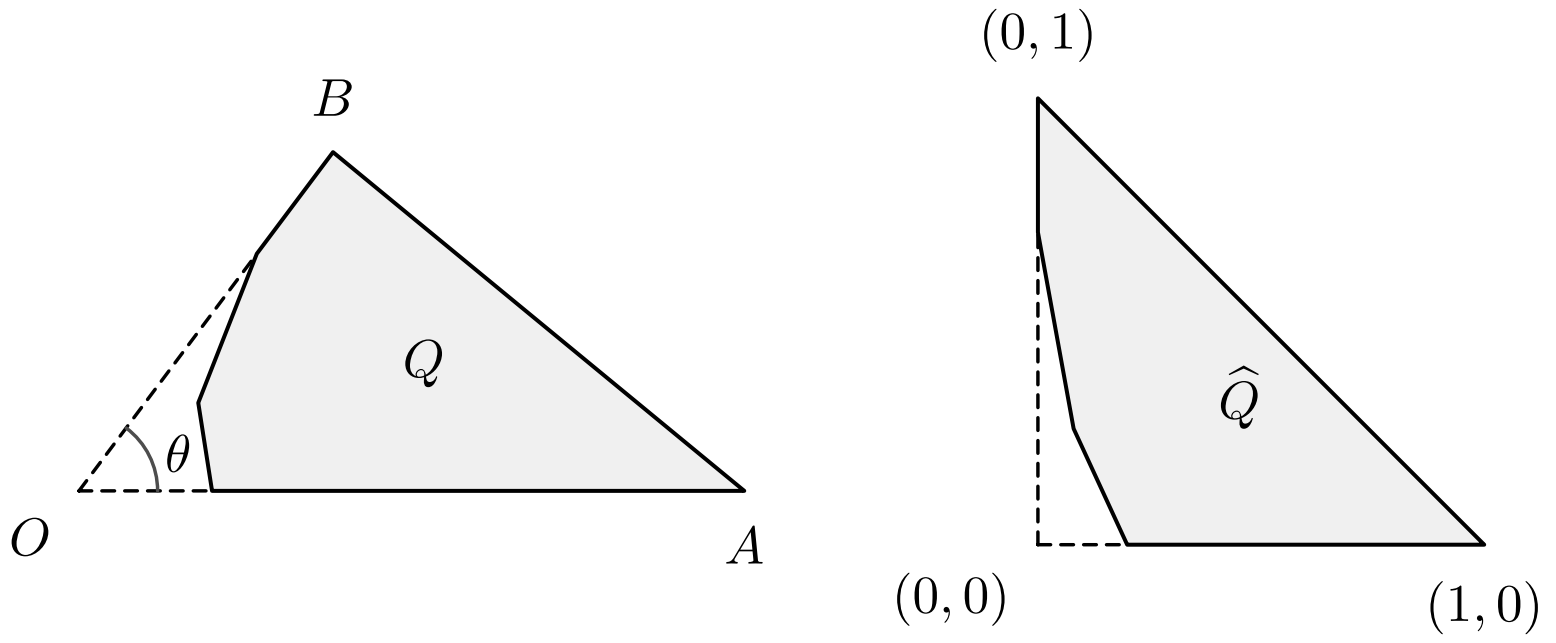
\par
Up to exchanging the two sides, we may assume $a\ge b$.
Moreover, we can scale the domain so that $ab=1$, which leaves
$\Rcal(Q)$ unchanged.

We observe that, given a feasible flux $F$ on $\widehat Q$, the field
$\sigma$ defined on $Q$ by
\begin{equation}\label{cc:eq:flux-pushforward}
 \sigma(X(x,y))=DX\,F(x,y),
\end{equation}
is a feasible flux on $Q$. Indeed, since $DX$ is constant,
$(\operatorname{div}\sigma)\circ X=\operatorname{div}F=-1$, and the
homogeneous Neumann condition is preserved.
A change of variables gives
\[
 \int_Q|\sigma|^2=\sin\theta\int_{\widehat Q}F^TH_{t,z}F,
\]
where $H_{t,z}$ is the matrix
\[
 H_{t,z}:=(DX)^TDX=\begin{pmatrix}t&z\\z&t^{-1}\end{pmatrix}.
\]
Then the upper bound in~\eqref{cc:eq:duality} yields
\begin{equation}\label{cc:eq:reference-upper}
 \Rcal(Q)\le\frac{4}{\delta^2\sin\theta}
 \int_{\widehat Q}F^TH_{t,z}F.
\end{equation}

\par\medskip\noindent-- \textit{Neumann-parallel flux bound.}
We now choose the field $F$ on $\widehat Q$ defined by
\[
 F(x,y)=(0,\ell(x)-y).
\]
Since only the second component of $F$ is nonzero, the field $\sigma$
associated with $F$ according to~\eqref{cc:eq:flux-pushforward} is parallel
to the shorter Neumann side $[O,B]$.

The field $F$ is feasible on $\widehat Q$. Indeed,
$\operatorname{div}F=-1$. The Neumann boundary consists of the lower
graph $y=\ell(x)$ and the portion on the $y$-axis. On the first,
$F=0$; on the second, $F$ is tangent to the boundary, since its first
component vanishes. Thus $F\cdot\nu=0$ on the whole Neumann boundary.
We have
\[
 \begin{aligned}
 \int_{\widehat Q}F^TH_{t,z}F
 &=\frac1t\int_0^1\int_{\ell(x)}^{1-x}(y-\ell(x))^2\dd y\dd x\\
 &=\frac1{3t}\int_0^1h(x)^3\dd x 
 \leq  \frac {2}{3t}   \left(\int_0^1 h\dd x\right)^3  =  \frac{\delta^3}{12t}\,,
  \end{aligned}
\]
where we have applied to $h$ the Favard--Berwald inequality for a
nonnegative concave function on $[0,1]$ (see~\cite{HeMa95}). 
Using~\eqref{cc:eq:reference-upper}, this yields~\eqref{cc:eq:column}.

\par\medskip\noindent-- \textit{Dirichlet-orthogonal flux bound.}
Assume that the orthogonal projection of $O$ onto the line through
$A$ and $B$ lies in the interior of $[A,B]$. We now normalize
$|AB|=1$, which leaves the ratio $t=a/b$ unchanged. After a rigid
motion, we may take
\[
 A=(0,0),\quad B=(1,0),\quad O=(\alpha,d),\qquad 0<\alpha<1.
\]
We use Cartesian coordinates $(r,s)$, where $r$ is measured along
$[A,B]$ and $s$ is the height above it. Then
$Q=\{(r,s):0<r<1,\ 0<s<\eta(r)\}$, where $\eta$ is nonnegative
and concave, and its graph is the Neumann boundary of $Q$.
Consider
\[
 \sigma(r,s):=(0,\eta(r)-s).
\]
It is a feasible flux, since $\operatorname{div}\sigma=-1$ in $Q$
and $\sigma=0$ on the Neumann boundary, where $s=\eta(r)$.
Applying the upper bound in~\eqref{cc:eq:duality} and integrating over
the vertical sections of $Q$, we obtain
\[
 \begin{aligned}
 \T_{\rm mix}(Q)
 &\le \int_Q|\sigma|^2\dd r\dd s
 =\int_0^1\int_0^{\eta(r)}(\eta(r)-s)^2\dd s\dd r\\
 &=\frac13\int_0^1\eta(r)^3\dd r
 =\int_0^1\int_0^{\eta(r)}s^2\dd s\dd r
 =\int_Q s^2\dd r\dd s.
 \end{aligned}
\]
Applying again the Favard--Berwald inequality  to $\eta$ and using $|Q|=\delta d/2$, we obtain
\[
 \T_{\rm mix}(Q)\le\frac23|Q|^3
 =\frac{\delta^3d^3}{12},\qquad
 \Rcal(Q)\le\frac{\delta d}{3}.
\]

It remains to express $d$ in terms of $\theta$ and $t$.
Before normalizing the length of $[A,B]$ to $1$, the distance from
$O$ to the line through $A$ and $B$ was $ab\sin\theta/|AB|$.
Thus the normalized height $d$ is $ab\sin\theta/|AB|^2$.
By the law of cosines, $|AB|^2=a^2+b^2-2abz$, so
\[
 d=\frac{ab\sin\theta}{a^2+b^2-2abz}
  =\frac{\sin\theta}{t+t^{-1}-2z}.
\]
Substituting this expression into $\Rcal(Q)\le\delta d/3$ gives
\[
 \Rcal(Q)\le\frac{\delta\sin\theta}{3(t+t^{-1}-2z)},
\]
which proves~\eqref{cc:eq:normal-roof}. The endpoint cases follow
by continuity.

For later use, we note that the condition on the orthogonal
projection of $O$ is equivalent to
\begin{equation}\label{cc:eq:projection-condition}
 1-\frac zt>0,\qquad 1-tz>0.
\end{equation}
Indeed, the first inequality ensures that the projection lies on
the $B$-side of $A$, while the second ensures that it lies on the
$A$-side of $B$. Since $t\ge1$ and $z<1$, the first inequality is
automatically satisfied. Thus, the condition reduces to $tz<1$.
In particular, it always holds when $z\le0$, whereas for $z>0$
it is equivalent to $t<1/z$.

\par\medskip\noindent-- \textit{Combining the two bounds.}
If $O$ projects orthogonally inside $[A,B]$, we can use
\eqref{cc:eq:normal-roof}; since $t+t^{-1}\geq 2$ we obtain
\[
 \Rcal(Q)\le\frac{\delta\sin\theta}{6(1-z)}
 =\frac\delta6\cot\frac\theta2.
\]
Otherwise $z>0$ and $t\ge1/z$, and
\eqref{cc:eq:column} yields
\[
 \Rcal(Q)\le\frac{\delta z}{3\sin\theta}
 \le\frac{\delta(1+z)}{6\sin\theta}
 =\frac\delta6\cot\frac\theta2.
\]

\end{proof}

\begin{proposition}[Lower bounds for $\Rcal(\theta)$]\label{cc:prop:large}
\begin{align}
 \Rcal(\theta)&\ge\frac18\cot\frac\theta2,\label{cc:eq:lower-benchmarks}\\[\smallskipamount]
 \Rcal(\theta)&\ge\frac{40\sin\theta}{9(32-27\cos\theta)}.\label{cc:eq:lower-benchmarks-second}
\end{align}
\end{proposition}
\begin{proof}
By scale invariance, we may assume that the two equal sides of
$I_\theta$ have unit length. The ray-coordinate map
$X(x,y)=xe_1+ye_2$
sends $\Delta$ onto $I_\theta$, with
\[
H_{1,z}= (DX) ^ T DX = (e_i\cdot e_j)_{i,j=1,2}
=\begin{pmatrix}
1 & z\\
z & 1
\end{pmatrix}\,, 
\qquad  \text{ with } z=\cos\theta.
\]
For any nonzero $v\in H^1(\Delta)$ vanishing on $x+y=1$,
we apply the variational lower bound in
\eqref{cc:eq:duality} to $v\circ X^{-1}$.
Since
\[
|\det DX|=\sin\theta,
\qquad |I_\theta|=\frac{\sin\theta}{2},
\]
a change of variable yields 
\begin{equation}\label{cc:eq:reference-trial}
\Rcal(\theta)\geq
\frac{4}{\sin\theta}
\frac{\left(\displaystyle\int_\Delta v\right)^2}
{\displaystyle\int_\Delta
\nabla v^T H_{1,z}^{-1}\nabla v}.
\end{equation}
Direct computations show that the lower bounds 
\eqref{cc:eq:lower-benchmarks} and \eqref{cc:eq:lower-benchmarks-second} 
 follow from  \eqref{cc:eq:reference-trial}  respectively taking as trial function  $v$ 
 the quadratic function 
 $$v_1(x, y)=1-(x+y)^2 $$ 
 and the quartic function
 $$
 v_2(x,y)=[1-(x+y)^2][1-(x-y)^2] \,. 
$$  
\end{proof}

\begin{corollary}[The large-side-ratio regime]\label{cc:cor:large-side-ratio}
If $t\ge2$, every $Q\in\mathcal C_\theta(T)$ satisfies $\Rcal(Q)<\Rcal(\theta)$.
\end{corollary}

\begin{proof}
We distinguish two cases.

\par\medskip\noindent-- \textit{Case 1: $-1<z<1/2$.}
At $t=2$, condition~\eqref{cc:eq:projection-condition} holds, so
inequality~\eqref{cc:eq:normal-roof} in Proposition~\ref{cc:prop:bounds} gives
\[
 \Rcal(Q)\le U(z):=\frac{2\sin\theta}{3(5-4z)}.
\]
Let us show that the same bound holds also for every $t\ge2$.
If $z\le0$, inequality~\eqref{cc:eq:normal-roof} in
Proposition~\ref{cc:prop:bounds} applies for all $t\ge2$,
and its right-hand side decreases with $t$.
If $0<z<1/2$, it applies for $2\le t<1/z$; for $t\ge1/z$,
we use inequality~\eqref{cc:eq:column} in
Proposition~\ref{cc:prop:bounds} instead.
The right-hand sides of the two estimates agree at $t=1/z$
and decrease with $t$ on their respective ranges.
Thus $\Rcal(Q)\le U(z)$. On the other hand,
inequality~\eqref{cc:eq:lower-benchmarks-second} in
Proposition~\ref{cc:prop:large} yields
\[
 \Rcal(\theta)-U(z)
 \ge\frac{40\sin\theta}{9(32-27z)}-U(z)
 =\frac{2\sin\theta(4+z)}{9(32-27z)(5-4z)}>0.
\]

\par\medskip\noindent-- \textit{Case 2: $1/2\le z<1$.}
For $t\ge2$, inequality~\eqref{cc:eq:column} in
Proposition~\ref{cc:prop:bounds} and inequality~\eqref{cc:eq:lower-benchmarks}
in Proposition~\ref{cc:prop:large} give
\[
 \Rcal(Q)\le\frac{\delta}{6\sin\theta}
 \le\frac1{6\sin\theta}
 <\frac{1+z}{8\sin\theta}\le\Rcal(\theta),
\]
where the strict inequality follows from $z\ge1/2$.
This completes the proof.
\end{proof}

\bigskip

\begin{corollary} \label{c:monobench}  The function $\theta\mapsto\Rcal(\theta)$ is strictly
 decreasing on $(0,\pi)$.
\end{corollary} 

\begin{proof}
Applying~\eqref{cc:eq:column} in Proposition~\ref{cc:prop:bounds} with $Q=T=I_\theta$ and recalling the definition of $\mathcal F$ in~\eqref{eq:triangular-F}, we obtain
\[
 0<\mathcal F(\theta)\le\frac{\theta^2}{3\sin\theta} \to 0 \qquad \text{ as } \theta \to 0 ^+\,.
\]
Hence, by Lemma \ref{cc:lem:triangle-benchmark}, $\mathcal F(\theta)/\theta$ is nonincreasing;
since it is positive,  the function
$$\Rcal(I_\theta)= \frac{ \mathcal F(\theta)} {\theta} \frac{1}{\theta}$$ is strictly decreasing.
\end{proof}

\subsection{\texorpdfstring{Polygonal first and second variations}{Polygonal first and second variations}}\label{cc:sec:variations}

For any convex polygon $Q$ with Dirichlet side $\Gamma = [A,B]$ (and Neumann
condition on the remaining sides), we 
consider for a fixed constant $c \in \R$ the functional
\[
 J(Q):=\T_{\rm mix}(Q)-c S \, , \qquad \text{ with } S: = |Q| \,. 
\]
We shall denote by  $\mathcal V$  the set of all vertices of $Q$ (including $A$ and $B$), and by $\mathcal E_N$
the set of Neumann edges.  

Our aim is to compute the first and second shape derivatives of $J$ under  a {\it linear vertex
velocity field $V$ vanishing on $\Gamma$}. By this, we mean that $V$ is a continuous piecewise
affine field on a triangulation of $Q$, which  agrees with a family of given velocity vectors
$V_Z$ at any $Z \in \mathcal V$  
and vanishes on $[A,B]$ (in particular, $V_A = V _B = 0$). Thus, under the action of $V$,  $[A,B]$ remains fixed, and every vertex different from $A, B$ moves linearly.  
We set 
\[
  Q_\varepsilon:=(\operatorname{Id}+\varepsilon V)(Q)\,.
 \]
The first and second order shape derivatives of $J$ at $Q$  are  defined respectively by
\begin{equation}\label{f:dder} 
 J'_Q[V]:= \left. \frac{d}{d\varepsilon}\right|_{\varepsilon=0} J ( Q _\varepsilon) \,, \qquad 
     J''_Q[V,V] :=  \left. \frac{d^2}{d\varepsilon^2}\right|_{\varepsilon=0} J ( Q _\varepsilon) \,.  
\end{equation}

We refer to~\cite{LNP16,LNP24,Laurain} for general results about
second-order shape derivatives on nonsmooth and convex domains; for elliptic energies with mixed
Dirichlet--Neumann conditions, including torsional rigidity,
see~\cite{BFL16} (in particular Example~3.10).
Although~\cite{BFL16} allow piecewise $C^1$ boundaries, including polygons,
their regularity assumptions on the solution and the deformation field
are not generally satisfied here: $w_Q$ need not be Lipschitz at
Dirichlet--Neumann vertices, and $V$ is only piecewise affine, rather than $C^1$.
We therefore give a short direct derivation for these vertex perturbations.
The first and second variation formulas are stated in Lemma~\ref{cc:lem:shape},
while the bulk-corner decomposition is established in Proposition~\ref{cc:prop:bulk-corner}.

\par

We denote by $\psi$ the {\it Eulerian shape derivative of $w_Q$},  defined  on  $Q$ by
\begin{equation}\label{f:dpsi}
\psi:=\dot w_Q-\nabla w_Q\cdot V, 
\end{equation}
 where $ \dot w_Q$ is  
the material derivative given on $Q$ by 
\[
\dot w _Q:=\left.\frac{d}{d\varepsilon}\right|_{\varepsilon=0}
(w_{Q_\varepsilon}\circ (\operatorname{Id}+\varepsilon V )) \,.
\]

We orient $\partial Q$ counterclockwise. On each $e\in\mathcal E_N$, let $t_e$
and $\nu_e$ be the oriented unit tangent and the outward unit normal. We set
\[
 V_t:=V\cdot t_e,\qquad V_\nu:=V\cdot\nu_e,
\]
and we write $(w _Q)_s$ and $(V_\nu)_s$  for the tangential derivative of $w_Q$  and $V _\nu$ along $e$;
notice that, by the Neumann condition, $|\nabla w_Q|^2=(w_Q)_s^2$ on every $e\in\mathcal E_N$.

The following lemma characterizes $\psi$ as the unique solution of a boundary value problem.
We point out that, when the additional regularity assumptions of ~\cite{BFL16}  hold,
this characterization allows to identify $\psi$ with the minimizer $v$ of the
nonlocal problem in~\cite[(3.11)]{BFL16}, specialized to mixed torsion: indeed, 
its Euler--Lagrange equation coincides with~\eqref{cc:eq:w-pde}.

\begin{lemma}[Boundary value problem satisfied by $\psi$]\label{cc:lem:eulerian}
For any convex polygon $Q$ and any linear vertex velocity field $V$
vanishing on $\Gamma=[A,B]$, the Eulerian derivative $\psi$
defined in~\eqref{f:dpsi} belongs to $H_\Gamma(Q)$ and is the unique weak solution of
\begin{equation}\label{cc:eq:w-pde}
 \begin{cases}
 -\Delta\psi=0 & \text{in }Q,\\
 \psi=0 & \text{on }\Gamma,\\
 \partial_\nu\psi=V_\nu+\bigl(V_\nu (w_Q)_s\bigr)_s
   & \text{on }\partial Q\setminus\Gamma.
 \end{cases}
\end{equation}
\end{lemma}

\begin{proof}

Set for brevity $w:=w_Q$ and $T_\varepsilon:=\operatorname{Id}+\varepsilon V$.
The pulled-back solution
$\widetilde w_\varepsilon=w_{Q_\varepsilon}\circ T_\varepsilon$
satisfies, on the fixed domain $Q$,
\begin{equation}\label{f:pulled} 
 \int_Q A_\varepsilon\nabla\widetilde w_\varepsilon\cdot
 \nabla\varphi\,\dd x
 =\int_Q\det(DT_\varepsilon)\varphi\,\dd x
 \qquad \forall \varphi\in H_\Gamma(Q),
\end{equation}
where
$A_\varepsilon:=\det(DT_\varepsilon)
(DT_\varepsilon)^{-1}(DT_\varepsilon)^{-T}$.
Since $V\in W^{1,\infty}(Q)$, both $A_\varepsilon$
and the right-hand side depend smoothly on $\varepsilon$.
Moreover, $A_\varepsilon\to I$ in $L^\infty(Q)$, so the
bilinear forms remain uniformly coercive on
$H_\Gamma(Q)$. 
By a standard perturbation argument for coercive
variational problems, the smooth dependence of the
bilinear form and the right-hand side on $\varepsilon$
implies that $\widetilde w_\varepsilon$ is
differentiable at $\varepsilon=0$ in $H_\Gamma(Q)$.

Let us check that the Eulerian derivative $\psi$ belongs
to the same space. By the corner expansions for mixed
elliptic problems~\cite[Chapter~4]{Grisvard}, the first
positive characteristic exponents of $w$ are
$\pi/\omega>1$ at Neumann--Neumann vertices and
$\pi/(2\omega)>1/2$ at Dirichlet--Neumann vertices,
where $\omega\in(0,\pi)$ is the interior angle.
Consequently, $D^2w\,V\in L^2(Q)$: near the former
vertices $D^2w\in L^2$, while near the latter
$V=O(r)$ compensates for the corner singularity. Since
\[
\nabla(\nabla w\cdot V)
=D^2w\,V+(DV)^T\nabla w
\]
and $V=0$ on $\Gamma$, it follows that
$\nabla w\cdot V\in H_\Gamma(Q)$. Hence
$\psi=\dot w-\nabla w\cdot V\in H_\Gamma(Q)$.

We now identify the equation satisfied by $\psi$.
Differentiating the variational identity~\eqref{f:pulled}
with respect to $\varepsilon$ at $0$, and then substituting
$\dot w=\psi+\nabla w\cdot V$ from~\eqref{f:dpsi}, we obtain
 \[
\int_Q\nabla\psi\cdot\nabla\varphi\,\dd x
=\int_Q(\operatorname{div}V)\varphi\,\dd x
 -\int_Q B(w,V)\cdot\nabla\varphi\,\dd x
 \qquad \forall \varphi\in H_\Gamma(Q),
\]
where
\[
B(w,V):=D^2w\,V+
 \bigl((\operatorname{div}V)I-DV\bigr)\nabla w.
\]

Since $D^2w\,V\in L^2(Q)$, as shown above, and
$DV\in L^\infty(Q)$, we have $B(w,V)\in L^2(Q)$.
 Moreover, $-\Delta w=1$ implies
$\operatorname{div}B(w,V)=-\operatorname{div}V$
in the distributional sense (also across the internal
edges of the triangulation, since the continuity of $V$ ensures
that the normal component of $B(w,V)$ has no jumps).

Integrating by parts the term involving $B(w,V)$,
the volume terms cancel by the divergence identity above, giving
\[
\int_Q\nabla\psi\cdot\nabla\varphi\,\dd x
=-\langle B(w,V)\cdot\nu,\varphi\rangle_{\partial Q}
 \qquad \forall \varphi\in H_\Gamma(Q).
\]

The weak normal trace of $B(w,V)$ is well defined since both
$B(w,V)$ and $\operatorname{div}B(w,V)$ belong to $L^2(Q)$.
 On each 
 Neumann edge, we have
$\Delta w=(w)_{ss}+(w)_{\nu_e\nu_e}$.
Moreover, the Neumann condition implies
$(w)_{s\nu_e}=0$.
 Hence,  
$$
-B(w,V)\cdot\nu_e
=-V_\nu (w)_{\nu_e\nu_e}+(w)_s(V_\nu)_s 
=V_\nu+\bigl(V_\nu(w)_s\bigr)_s.
$$
The corner estimates ensure that the identity, established on each open Neumann edge, also holds on the entire Neumann boundary in the weak sense. 
Uniqueness follows from the nonempty Dirichlet portion $\Gamma$.
\end{proof}

\begin{lemma}[Nontriviality of $\psi$]\label{cc:lem:strict-energy}
Let $Q$ be a convex polygon with Dirichlet side $\Gamma = [A,B]$.
If $V$ is a linear vertex velocity field vanishing on $\Gamma$
and $V_Z\ne0$ at some genuine vertex $Z\in\mathcal V\setminus\{A,B\}$,
then the Eulerian derivative $\psi$ defined in~\eqref{f:dpsi} is nonzero. In particular,
 \begin{equation}\label{cc:eq:strict-energy}
 \int_Q|\nabla\psi|^2\,\dd x>0.
\end{equation}
\end{lemma}

\begin{proof}
Assume by contradiction that $\int_Q|\nabla\psi|^2\,\dd x=0$.
Since $\psi=0$ on $\Gamma$, we have $\psi\equiv0$.
By~\eqref{cc:eq:w-pde}, on every $e\in\mathcal E_N$,
\begin{equation}\label{f:EN}
\bigl(V_\nu(w_Q)_s\bigr)_s=-V_\nu.
\end{equation} 
The corner estimates recalled in the proof of
Lemma~\ref{cc:lem:eulerian} ensure that
$V_\nu(w_Q)_s$ vanishes at both endpoints of $e$.
Indeed, at Dirichlet--Neumann vertices $V=O(r)$
compensates for the possible singularity of
$\nabla w_Q$, while at Neumann--Neumann vertices
$\nabla w_Q$ vanishes. Integrating along $e$,
we therefore obtain
\[
\int_e V_\nu\,\dd s
=-[V_\nu(w_Q)_s]_{\partial e}=0.
\]
If the affine function $V_\nu$ were not identically
zero, it would have a unique zero $s_0$ in the
interior of $e$. 
Integrating \eqref{f:EN}
from the initial endpoint of $e$ to $s_0$,
and using the vanishing of $V_\nu(w_Q)_s$
at both endpoints of this subsegment, we obtain
\[
\int_0^{s_0}V_\nu\,\dd s=0\,,
\]
which is impossible since $V_\nu$ has a constant
nonzero sign on such subsegment. Thus $V_\nu=0$
on every Neumann edge. At each
$Z\in\mathcal V\setminus\{A,B\}$,  the two incident
outward normals are linearly independent, so
$V_Z=0$, a contradiction.
\end{proof}

\bigskip 
On $\partial Q \setminus \Gamma$,  define
\begin{equation}\label{f:defvarie}
 f:=2w_Q- (w_Q)_s ^2-c\,.
\end{equation}
Moreover, if  $g$ is defined on an oriented edge $e\in \mathcal E _N$, with initial and terminal
vertices $Z_-$ and $Z_+$, set
$[g]_{\partial e}:=g(Z_+)-g(Z_-)$.

\begin{lemma}[Shape derivatives of $J$]\label{cc:lem:shape}
Let $Q$ be a convex polygon with Dirichlet side $\Gamma=[A,B]$,
and let $V$ be a linear vertex velocity field vanishing on $\Gamma$.
Then the first and second shape derivatives of $J$, defined in~\eqref{f:dder},
are given by

\begin{eqnarray}
&  \displaystyle J'_Q[V]=\int_{\partial Q\setminus\Gamma} fV_\nu\,\dd s  & \label{cc:eq:first} 
\\ \noalign{\medskip} 
 &  \displaystyle J''_Q[V,V] =2\int_Q|\nabla\psi|^2\,\dd x
  +\sum_{e\in\mathcal E_N}[fV_tV_\nu]_{\partial e}
  -2\int_{\partial Q\setminus\Gamma} fV_t(V_\nu)_s\,\dd s \,.  & \label{cc:eq:polygon-hessian}
 \end{eqnarray}
For the special choice $c=2\T_{\rm mix}(Q)/S$, we also have the following identities, where
all shape derivatives are computed at $Q$ along $V$: 
 \begin{equation}\label{cc:eq:normalization}
 S^2\Rcal'=J',\qquad
 S^2\Rcal''=J''-4\frac{S'}S J'-\frac cS(S')^2\,.
 \end{equation}
\end{lemma}

\begin{proof}
Set for brevity $w:=w_Q$ and $T_\varepsilon:=\operatorname{Id}+\varepsilon V$.
As in the proof of Lemma~\ref{cc:lem:eulerian}, the pulled-back
problem~\eqref{f:pulled} depends smoothly on $\varepsilon$ and remains
uniformly coercive. Hence $J(Q_\varepsilon)$ is twice differentiable
at $\varepsilon=0$.
Differentiating
$\T_{\rm mix}(Q_\varepsilon)
=\int_Q\widetilde w_\varepsilon
\det(DT_\varepsilon)\,\dd x$ at $\varepsilon=0$
gives
\[
\T_{\rm mix}'(Q)
=\int_Q\bigl(\dot w+w\operatorname{div}V\bigr)\,\dd x.
\]

To eliminate $\dot w$, we differentiate the
variational identity~\eqref{f:pulled} at $\varepsilon=0$
and choose $\varphi=w$. We then use the same
identity at $\varepsilon=0$ with $\varphi=\dot w$.
Using
$A'_0=(\operatorname{div}V)I-DV-DV^T$,
we obtain the following formula for the first variation of mixed torsional rigidity (which  agrees with~\cite[Example~3.10]{BFL16}, up to the factor $2$
in our normalization of  $\T_{\rm mix}$): 
\[  \T_{\rm mix}'(Q)
 =\int_Q\!\left[(2w-|\nabla w|^2)\operatorname{div}V
       +2\nabla w\cdot DV\,\nabla w\right] \,\dd x
 =\int_{\partial Q\setminus\Gamma}(2w-w_s^2)V_\nu\,\dd s\,.
\]
Here the second equality follows by applying the divergence theorem to
$(2w-|\nabla w|^2)V+2(V\cdot\nabla w)\nabla w$, using $V=0$
on $\Gamma$ and $\partial_\nu w=0$ on the Neumann boundary.

As $S'=\int_{\partial Q\setminus\Gamma}V_\nu\,\dd s$,
we obtain~\eqref{cc:eq:first}.

For the second variation we work entirely with the pulled-back bulk
energy; in particular, no differentiation of singular edge integrals
is assumed. Put
\[
a_\varepsilon(u,v)=\int_Q A_\varepsilon\nabla u\cdot\nabla v,
\qquad
\ell_\varepsilon(v)=\int_Q {\rm det} ( I + \varepsilon DV)  v.
\]
Writing $A_\varepsilon=I+\varepsilon A_1+\varepsilon^2A_2+O(\varepsilon^3)$,
we have, almost everywhere on each affine piece,
\[
\begin{aligned}
A_1&=(\operatorname{tr}DV)I-DV-DV^T,\\
A_2&=(\det DV)I-(\operatorname{tr}DV)(DV+DV^T)
       +DV^2+DVDV^T+(DV^T)^2.
\end{aligned}
\]
The coercive problem~\eqref{f:pulled} is $C^2$ in $H_\Gamma(Q)$.
Since $a_\varepsilon(\widetilde w_\varepsilon,v)
=\ell_\varepsilon(v)$, differentiation gives
\[
\T_{\rm mix}'=2\ell'_0(w)-a'_0(w,w),\qquad
\T_{\rm mix}''=2\ell''_0(w)-a''_0(w,w)
                    +2a_0(\dot w,\dot w).
\]
Consequently, as $|Q_\varepsilon|''|_{\varepsilon=0}
=2\int_Q\det DV$,
\begin{equation}\label{cc:eq:bulk-second}
J''_Q[V,V]
=2\int_Q|\nabla\dot w|^2
 +2\int_Q(2w-c)\det DV
 -2\int_Q\nabla w\cdot A_2\nabla w.
\end{equation}
The above identity involves only volume integrals over $Q$. 
To reduce it to boundary terms, substitute
$\dot w=\psi+V\cdot\nabla w$ and use the weak equation
\eqref{cc:eq:w-pde} with test function $V\cdot\nabla w\in
H_\Gamma(Q)$:
\[
\int_Q\nabla\psi\cdot\nabla(V\cdot\nabla w)
=\sum_{e\in\mathcal E_N}\int_e
 V_t w_s\bigl[V_\nu+(V_\nu w_s)_s\bigr].
\]
Expanding the remaining integrals in~\eqref{cc:eq:bulk-second},
using $-\Delta w=1$ and integrating by parts on each affine piece of
$V$ gives, after cancellation of the internal triangulation traces,
Indeed, across an internal edge with unit normal $n$, continuity of the piecewise affine field makes cancel the normal jumps of the first- and second-order fluxes arising in the expansion.
\[
\begin{aligned}
&2\int_Q|\nabla(V\cdot\nabla w)|^2
 +4\sum_e\int_e V_t w_s\bigl[V_\nu+(V_\nu w_s)_s\bigr]\\
&\quad+2\int_Q(2w-c)\det DV
 -2\int_Q\nabla w\cdot A_2\nabla w\\
&=\sum_{e\in\mathcal E_N}[fV_tV_\nu]_{\partial e}
 -2\sum_{e\in\mathcal E_N}\int_e fV_t(V_\nu)_s.
\end{aligned}
\]
Here the integrations are first performed with a radius-$r$ corner
neighbourhood removed. At Dirichlet--Neumann corners the exponent
$\beta>1/2$ and $V=O(r)$ give
$fV_tV_\nu=O(r^{2\beta}|\log r|^C)+O(r^2)\to0$;
at Neumann--Neumann corners $w_s\to0$, so
$f\to2w(Z)-c$. All edge integrands are locally integrable, and the
terms $V_\nu w_s\psi$ tend to zero by
Remark~\ref{cc:rem:psi-corners}. Thus the corner cutoffs can be removed.
Combining this identity with~\eqref{cc:eq:bulk-second} yields
\eqref{cc:eq:polygon-hessian}, without exchanging differentiation
and integration on a singular edge.
\par

Finally, for $c=2\T_{\rm mix}(Q)/S$ held fixed,
differentiating $\Rcal=\T_{\rm mix}/S^2$ twice and using
$\T_{\rm mix}'=J'+cS'$ and $\T_{\rm mix}''=J''+cS''$
yields~\eqref{cc:eq:normalization}.
\end{proof}

\begin{remark}[Corner estimates for $\psi$]\label{cc:rem:psi-corners}
Set $g:=V_\nu+(V_\nu(w_Q)_s)_s$, the Neumann datum in~\eqref{cc:eq:w-pde}.
At a Dirichlet--Neumann vertex of angle $\omega$, let
$\beta:=\min\{2,\pi/(2\omega)\}>1/2$. Standard corner expansions
(see~\cite[Chapter~4]{Grisvard}) give
\[
(w_Q)_s=O(r^{\beta-1}|\log r|^k),\qquad
g=O(r^{\beta-1}|\log r|^k)
\]
for some $k\ge0$. If $\max\{0,1-\beta\}<\gamma<\beta$, the corresponding
local mixed-boundary estimates yield
$\psi=O(r^\gamma)$ and $\psi_s=O(r^{\gamma-1})$. Hence
$(w_Q)_s\psi\to0$ and $V_\nu(w_Q)_s\psi_s\in L^1$.

At a Neumann--Neumann vertex, with
$\alpha:=\min\{2,\pi/\omega\}>1$, one has
$(w_Q)_s=O(r^{\alpha-1}|\log r|^k)$ and
$g=O(r^{\alpha-2}|\log r|^k)$. The analogous local estimates give
$\psi=O(1)$ and $\psi_s=O(r^{-1})$, so the same two conclusions hold.
Therefore all integrations by parts in Lemma~\ref{cc:lem:shape} are justified and
$[V_\nu(w_Q)_s\psi]_{\partial e}=0$ on every Neumann edge.
\end{remark}

\begin{proposition}[Bulk-corner decomposition of the second variation at stationary polygons]\label{cc:prop:bulk-corner}
Let $Q$ be a convex polygon with Dirichlet side $\Gamma=[A,B]$.
Let $\mathcal M\subset\mathcal V\setminus\{A,B\}$, and assume that
$J'_Q[V]=0$ whenever a single vertex $Z\in\mathcal M$ moves with
arbitrary velocity $V_Z\in\R^2$, all other vertices being fixed.
If $V$ is a linear vertex velocity field vanishing on $\Gamma$
and at all vertices outside $\mathcal M$, then
\begin{equation}\label{cc:eq:corner-hessian}
 J''_Q[V,V]=2\int_Q|\nabla\psi|^2\,\dd x
 +\sum_{Z\in\mathcal M}(2w_Q(Z)-c)\mathcal Q_Z(V_Z).
\end{equation}
 For each $Z\in\mathcal M$, 
the quadratic form $\mathcal Q_Z$ is represented in the orthonormal
basis $(t_{e_-(Z)},-\nu_{e_-(Z)})$  (where 
$e_-(Z)$ denotes the Neumann edge
entering $Z$) by the matrix 

\begin{equation}\label{cc:eq:corner-matrix}
 M_Z=\sin\omega_Z
 \begin{pmatrix}\cos\omega_Z&-\sin\omega_Z\\
                 -\sin\omega_Z&-\cos\omega_Z\end{pmatrix},
\end{equation}
 where $\omega_Z\in(0,\pi)$ is the interior angle of $Q$ at $Z$. In particular, $\det M_Z=-\sin^2\omega_Z<0$.
\end{proposition}

\begin{proof}

For each vertex $Z$ of $Q$, let $\eta_Z$ be the
boundary hat function, affine on each edge, equal
to $1$ at $Z$ and $0$ at all other vertices.

Fix $Z\in\mathcal M$, and abbreviate
$t_\pm=t_{e_\pm(Z)}$ and $\nu_\pm=\nu_{e_\pm(Z)}$.
When only $Z$ moves, the first-variation
formula~\eqref{cc:eq:first} gives
\[
J'_Q[V]=V_Z\cdot\left(
\nu_-\int_{e_-(Z)}f\eta_Z\,\dd s+
\nu_+\int_{e_+(Z)}f\eta_Z\,\dd s
\right).
\]
Since $\nu_-$ and $\nu_+$ are linearly independent,
stationarity with respect to arbitrary $V_Z$ implies
\begin{equation}\label{cc:eq:zero-hats}
\int_{e_-(Z)}f\eta_Z\,\dd s
=\int_{e_+(Z)}f\eta_Z\,\dd s=0.
\end{equation}

For a general velocity field supported on $\mathcal M$,
$V_t$ is a linear combination of the corresponding
hat functions on each edge $e$. Thus
$\int_e fV_t\,\dd s=0$ by~\eqref{cc:eq:zero-hats}.
Since $(V_\nu)_s$ is constant along each edge,
the integral term in~\eqref{cc:eq:polygon-hessian}
vanishes:
\[
-2\int_{\partial Q\setminus\Gamma}
fV_t(V_\nu)_s\,\dd s
=-2\sum_{e\in\mathcal E_N}(V_\nu)_s
\int_e fV_t\,\dd s=0.
\]

It remains to evaluate the endpoint term
\[
\sum_{e\in\mathcal E_N}
[fV_tV_\nu]_{\partial e}.
\]
Since the edges are oriented counterclockwise,
$e_-(Z)$ ends at $Z$, whereas $e_+(Z)$ starts
at $Z$. Thus the contributions of these two
edges have opposite signs, and their sum is
\[
f(Z)\bigl[
(V_Z\cdot t_-)(V_Z\cdot\nu_-)
-(V_Z\cdot t_+)(V_Z\cdot\nu_+)
\bigr].
\]
All other vertex contributions vanish because
$V$ is supported on $\mathcal M$.
Moreover, both edges incident to $Z$ are Neumann,
and convex-corner regularity gives
$\nabla w_Q(Z)=0$, hence $f(Z)=2w_Q(Z)-c$.
This proves~\eqref{cc:eq:corner-hessian}.

To compute the matrix of $\mathcal Q_Z$, take
coordinates $(x,y)$ in the orthonormal basis
$(t_-,-\nu_-)$ and set $\omega=\omega_Z$.
Then
$t_+=(-\cos\omega,\sin\omega)$ and 
$\nu_+=(\sin\omega,\cos\omega)$, so that 
\[
\begin{aligned}
\mathcal Q_Z(x,y)
&=-xy-(x\sin\omega+y\cos\omega)
             (-x\cos\omega+y\sin\omega)\\
&=\sin\omega\bigl(
\cos\omega(x^2-y^2)-2\sin\omega\,xy
\bigr).
\end{aligned}
\]
This yields the matrix~\eqref{cc:eq:corner-matrix},
whose determinant is
$-\sin^2\omega_Z<0$.
\end{proof}

\subsection{\texorpdfstring{{Reduction to five-vertex polygonal caps}}{Reduction to five-vertex polygonal caps}}\label{cc:sec:five}
\leavevmode\par

Given a triangle $T=\triangle(O,A,B)$, let $\mathcal P_\theta(T)$  be the sub-class of $ \mathcal C _\theta (T)$ consisting of  five-vertex polygonal caps

\begin{equation}\label{cc:eq:five-vertex-class}
\begin{aligned}
\mathcal P_\theta(T):=\bigl\{P:\ &\overline P=\conv\{A,C,X,D,B\},\\
&C\in[O,A],\quad D\in[O,B],\\
&X\in\overline T,\quad |P|>0\bigr\}.
\end{aligned}
\end{equation}

\begin{lemma}[Continuity of mixed torsion under cap convergence]\label{cc:lem:cap-continuity}
Let $\Gamma$ be a fixed nondegenerate segment and let $Q_k$ be uniformly bounded
convex domains lying on the same side of the line containing $\Gamma$, with
$\Gamma\subset\partial Q_k$. If $\overline Q_k$ converges in the Hausdorff metric
to the closure of a convex domain $Q$ of positive area and
$\Gamma\subset\partial Q$, then
\[
\T_{\rm mix}(Q_k,\Gamma)\longrightarrow\T_{\rm mix}(Q,\Gamma).
\]
\end{lemma}
\begin{proof}
Embed the domains in a fixed bounded Lipschitz set. Convexity, uniform boundedness
and the positive area of the limit give, for large $k$, a uniform interior cone
condition and hence uniformly bounded extension operators; see
\cite[Section~3.7]{HePi18}. The fixed segment $\Gamma$ gives a uniform Poincar\'e
inequality, and weak limits preserve the zero trace on $\Gamma$.
Since convex Hausdorff convergence implies convergence of the characteristic
functions almost everywhere, weak compactness and lower semicontinuity give one
energy inequality. For the converse, approximate functions in $H_\Gamma(Q)$ by
restrictions of smooth functions vanishing on $\Gamma$, restrict them to $Q_k$,
and use a diagonal argument. The variational characterization~\eqref{eq:mixed-primal}
then yields the claim.
\end{proof}

\begin{theorem}[Reduction to five-vertex polygonal caps]\label{cc:thm:five}
For every $Q\in\mathcal C_\theta(T)$, there exists  $P\in\mathcal P_\theta(T)$
such that
\[
 |P|=|Q|,\qquad \T_{\rm mix}(P)\ge\T_{\rm mix}(Q).
\]
Consequently,
\begin{equation}\label{cc:eq:five}
 \sup\bigl\{\Rcal(Q):Q\in\mathcal C_\theta(T)\bigr\}
 =\max\bigl\{\Rcal(P):P\in\mathcal P_\theta(T)\bigr\}.
\end{equation}
\end{theorem}
\begin{proof}
First, we observe that the maximum on the right of~\eqref{cc:eq:five} is attained.
Indeed, $\Rcal$ is continuous on positive-area caps by
Lemma~\ref{cc:lem:cap-continuity}, and it extends continuously by zero to zero-area caps by
Proposition~\ref{cc:prop:bounds}, since $\theta$ is fixed.
The family of possibly degenerate five-vertex caps is compact; hence the
maximum is attained at a cap of positive area.
We now fix $s\in(0,|T|]$: within the class of caps of area $s$, maximizing $\Rcal$ is equivalent to maximizing $\T_{\rm mix}$.

Let $Q\in\mathcal C_\theta(T)$ have area $s$, and choose $N$ not
smaller than the number of vertices of $Q$.

We observe that the problem of maximizing the mixed torsional rigidity among polygonal caps in $\mathcal C _\theta (T)$
having area $s$ and at most $N$ vertices admits a solution $Q_N$.
Indeed, this class is compact under Hausdorff convergence,
and the positive fixed area ensures the continuity of the mixed torsional rigidity. 

Then we are reduced to show that $Q_N\in\mathcal P_\theta(T)$.
To that aim, it is enough to prove that at most one
of its genuine vertices lies in $T$.
Suppose, by contradiction, that $Q_N$ has two genuine
vertices $Z_1,Z_2\in T$. 
Since these are genuine vertices lying in $T$,
sufficiently small independent displacements in any direction
preserve admissibility. 
Moreover, the first variation of the area with respect
to either vertex is not identically zero. Hence the
area constraint  is regular, and the Lagrange
multiplier rule yields $c\in\mathbb R$ such that
\[
J=\mathcal T_{\rm mix}-cS
\]
is stationary at $Q_N$ under independent variations of
$Z_1$ and $Z_2$.

By Proposition~\ref{cc:prop:bulk-corner}, the quadratic
forms $\mathcal Q_{Z_1}$ and $\mathcal Q_{Z_2}$ are indefinite.
We can therefore choose nonzero vectors $E_1,E_2\in\mathbb R^2$
such that
\[
(2w_{Q_N}(Z_i)-c)\mathcal Q_{Z_i}(E_i)\geq 0,
\qquad i=1,2.
\]

We now consider a velocity field $V$ that
leaves all other vertices
fixed and 
moves $Z_i$
with velocity $\alpha_i E_i$, where the coefficients
$\alpha _i$ are chosen so that  $S'_{Q_N}[V]=0$. (Notice that the latter  is a single homogeneous
linear equation in $\alpha_1,\alpha_2$,  so that we can find a solution
with  not both $\alpha _i$ equal to zero).  
For this choice of $V$,  using
Proposition~\ref{cc:prop:bulk-corner} and
Lemma~\ref{cc:lem:strict-energy}, we obtain 
 \[
\begin{aligned}
J''_{Q_N}[V,V]
&=2\int_{Q_N}|\nabla\psi|^2\,dx\\
&\quad+\sum_{i=1}^{2}
(2w_{Q_N}(Z_i)-c)\,
\mathcal Q_{Z_i}(\alpha_i E_i)>0.
\end{aligned}
\]

Choose $W$ supported at $Z_1$ with $S'_{Q_N}[W]\ne0$,
and let $Q_{\varepsilon,\delta}$ be the polygon obtained by
moving the vertices of $Q_N$ by $\varepsilon V+\delta W$.

Apply the implicit function theorem to
$F(\varepsilon,\delta):=|Q_{\varepsilon,\delta}|-s$.
Since $F(0,0)=0$, $\partial_\delta F(0,0)=S'_{Q_N}[W]\ne0$
and $\partial_\varepsilon F(0,0)=S'_{Q_N}[V]=0$,
there exists a smooth $\eta(\varepsilon)$ such that
$|Q_{\varepsilon,\eta(\varepsilon)}|=s$ and
$\eta(0)=0$.
Differentiating the area constraint
$F(\varepsilon,\eta(\varepsilon))=0$ at $\varepsilon=0$ gives
$0=S'_{Q_N}[V]+\eta'(0)S'_{Q_N}[W]$. Since
$S'_{Q_N}[V]=0$ and $S'_{Q_N}[W]\ne0$, we obtain
$\eta'(0)=0$.
Along this area-preserving deformation,  applying the chain rule 
and using the stationarity of $J$ at $Q_N$,  we obtain 
\[
\begin{aligned}
\left.\frac{\dd^2}{\dd\varepsilon^2}
\T_{\rm mix}(Q_{\varepsilon,\eta(\varepsilon)})
\right|_{\varepsilon=0}
&=J''_{Q_N}[V,V]+\eta''(0)J'_{Q_N}[W]\\
&=J''_{Q_N}[V,V]>0.
\end{aligned}
\]

Moreover, since $\eta'(0)=0$, the corrected perturbation
has initial velocity $V$. Hence
\[
\left.\frac{d}{d\varepsilon}
\mathcal T_{\rm mix}(Q_{\varepsilon,\eta(\varepsilon)})
\right|_{\varepsilon=0}
=J'_{Q_N}[V]+cS'_{Q_N}[V]=0,
\]
by the stationarity of $J$ at $Q_N$ and the choice
of $V$. Thus the mixed torsion increases for small
nonzero $\varepsilon$, contradicting the maximality
of $Q_N$.

Hence $Q_N\in\mathcal P_\theta(T)$. Since $Q$ is admissible for the
$N$-vertex maximization problem,
\[
\T_{\rm mix}(Q)\leq\T_{\rm mix}(Q_N),
\qquad |Q|=|Q_N|.
\]

\end{proof}

\subsection{\texorpdfstring{Reduction to triangular caps}{Reduction to triangular caps}}\label{cc:sec:exclusion}\label{cc:sec:angular-motion}
Fix a reference triangle $T=\triangle(O,A,B)$ with opening angle
$\theta=\angle AOB\in(0,\pi)$ and consider a five-vertex polygonal cap
$\overline P=\conv\{A,C,X,D,B\}$ for $P\in\mathcal P_\theta(T)$, with
$\mathcal P_\theta(T)$ defined by~\eqref{cc:eq:five-vertex-class}.
Introduce the weighted functional
\begin{equation}\label{cc:eq:Phi}
 \Phi(P,T)=\theta^2\Rcal(P)-\mathcal F(\theta)
          =\theta^2\bigl(\Rcal(P)-\Rcal(\theta)\bigr).
\end{equation}
We study its second derivative when  $A$ and $B$ are kept fixed but, besides the remaining vertices, also  the apex angle $\theta$ of the containing
triangle $T$ is allowed to vary.

\begin{figure}[H]
\centering
\begin{minipage}[t]{.32\linewidth}\centering
\begin{tikzpicture}[x=.60cm,y=.60cm,font=\footnotesize,line join=round,line cap=round,draw=black,every node/.append style={text=black}]
\coordinate (O) at (0,4); \coordinate (A) at (-2.7,0); \coordinate (B) at (2.7,0);
\coordinate (C) at (-1.6,1.6296); \coordinate (D) at (1.6,1.6296); \coordinate (X) at (0,2.1);
\draw[psDashed] (A)--(O)--(B);
\filldraw[psHull,psOutline] (A)--(C)--(X)--(D)--(B)--cycle;
\draw[psInternal] (O) ++(-124:.48) arc[start angle=-124,end angle=-56,radius=.48];
\node[above=2pt] at (O) {$O$}; \node[below left=1pt] at (A) {$A$}; \node[below right=1pt] at (B) {$B$};
\node[left=1pt] at (C) {$C$}; \node[right=1pt] at (D) {$D$}; \node[above=1pt] at (X) {$X$};
\node at (0,3.20) {$\theta$};
\end{tikzpicture}
\end{minipage}\hfill
\begin{minipage}[t]{.32\linewidth}\centering
\begin{tikzpicture}[x=.60cm,y=.60cm,font=\footnotesize,line join=round,line cap=round,draw=black,every node/.append style={text=black}]
\coordinate (O) at (0,4); \coordinate (A) at (-2.7,0); \coordinate (B) at (2.7,0);
\coordinate (C) at (-1.6,1.6296); \coordinate (D) at (1.85,1.2593);
\draw[psDashed] (A)--(O)--(B);
\filldraw[psHull,psOutline] (A)--(C)--(D)--(B)--cycle;
\draw[psInternal] (O) ++(-124:.48) arc[start angle=-124,end angle=-56,radius=.48];
\node[above=2pt] at (O) {$O$}; \node[below left=1pt] at (A) {$A$}; \node[below right=1pt] at (B) {$B$};
\node[left=1pt] at (C) {$C$}; \node[right=1pt] at (D) {$D$};
\node at (0,3.20) {$\theta$};
\end{tikzpicture}
\end{minipage}\hfill
\begin{minipage}[t]{.32\linewidth}\centering
\begin{tikzpicture}[x=.60cm,y=.60cm,font=\footnotesize,line join=round,line cap=round,draw=black,every node/.append style={text=black}]
\coordinate (O) at (0,4); \coordinate (A) at (-2.7,0); \coordinate (B) at (2.7,0);
\coordinate (X) at (-.6,2.1); \coordinate (D) at (1.6,1.6296);
\draw[psDashed] (A)--(O)--(B);
\filldraw[psHull,psOutline] (A)--(X)--(D)--(B)--cycle;
\draw[psInternal] (O) ++(-124:.48) arc[start angle=-124,end angle=-56,radius=.48];
\node[above=2pt] at (O) {$O$}; \node[below left=1pt] at (A) {$A$}; \node[below right=1pt] at (B) {$B$};
\node[above left=2pt] at (X) {$X$}; \node[right=1pt] at (D) {$D$};
\node at (0,3.20) {$\theta$};
\end{tikzpicture}
\end{minipage}
\caption{The three possible nontriangular configurations of a five-vertex
polygonal cap.}
\label{cc:fig:three-genuine-caps}
\end{figure}
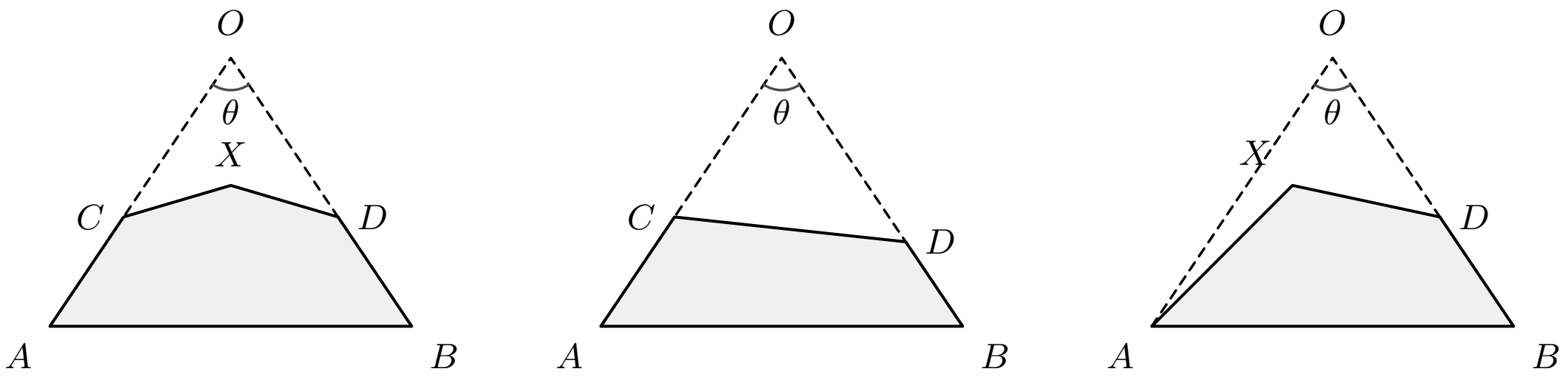

Unless $P$ is a triangular cap, we can distinguish the following three
configurations (see Figure~\ref{cc:fig:three-genuine-caps}):
 \begin{enumerate}[label=(\roman*),leftmargin=*,itemsep=2pt]
\item a \textit{pentagonal cap}
$\overline P=\conv\{A,C,X,D,B\}$, with $C\in(A,O)$, $D\in(B,O)$ and
$X\in T$;
\item a \textit{two-contact quadrilateral cap}
$\overline P=\conv\{A,C,D,B\}$, with $C\in(A,O)$ and $D\in(B,O)$;
\item a \textit{one-contact quadrilateral cap}
$\overline P=\conv\{A,X,D,B\}$, with $X\in T$, $D\in(B,O)$ and
$[O,A]\cap\overline P=\{A\}$, or its reflection with $A$ and $B$ interchanged.
\end{enumerate}

In each of the above cases, the {\it free Euclidean coordinates} will be: 
\smallskip 

\begin{itemize} 

\item[--] In Case~(i), the Cartesian coordinates of the free vertices
$(C,X,D)\in\mathbb R^6$.
\smallskip 

\item[--]  In Case~(ii), the Cartesian coordinates of the free vertices 
$(C,D)\in\mathbb R^4$.

\smallskip 
\item[--]  In Case~(iii), the Cartesian coordinates of the free vertices $(X, D)$, together
with the opening angle,  so $(X,D,\theta)\in\mathbb R^5$.

\end{itemize} 
We study the second variation of $\Phi$ with respect to these free Euclidean coordinates.
The following proposition covers all three cases. In Cases~(i) and~(ii),
we allow $\theta$ to vary; in Case~(iii), we restrict to variations of
$(X,D)$ with $\theta$ fixed. For a linear vertex velocity field $V$
vanishing on $[A,B]$ (as defined in Section~\ref{cc:sec:variations}),
we write 
\[
 \Phi''_{(P,T)}[V,V]
 :=\left.\frac{d^2}{d\varepsilon^2}\right|_{\varepsilon=0}
 \Phi(P_\varepsilon,T_\varepsilon)
\]
Here $P_\varepsilon=(\operatorname{Id}+\varepsilon V)(P)$.
In Cases~(i) and~(ii), $T_\varepsilon$ is the containing triangle
whose sides lie on the lines $AC_\varepsilon$ and $BD_\varepsilon$,
where $C_\varepsilon=C+\varepsilon V_C$ and
$D_\varepsilon=D+\varepsilon V_D$. In Case~(iii),
$T_\varepsilon$ is the containing triangle determined locally by
$A,B,D_\varepsilon$ and the prescribed angle $\theta$.

\par\smallskip

\begin{proposition}[Structure of the Hessian at stationary configurations]
\label{cc:prop:stationary-hessian}
Let $P\in\mathcal P_\theta(T)$ be a nontriangular cap in one of the
three configurations (i)--(iii) described above, and set $S=|P|$.
Assume that $(P,T)$ is stationary for $\Phi$ with respect to the
free Euclidean coordinates, with $\theta$ held fixed in Case~(iii).
Then, for every corresponding linear vertex velocity field $V$
vanishing on $[A,B]$, 
\begin{equation}\label{cc:eq:weighted-hessian}
 S^2\Phi''_{(P,T)}[V,V]
 =2\theta^2\int_P|\nabla\psi|^2\,\dd x+q_0(V)+q_r(V)\,.
\end{equation}
Here, denoting by $\mathcal M$ the set
of free vertices, the quadratic form $q_0$ is of the form
\begin{equation}\label{f:q0} 
 q_0(V)=\sum_{Z\in\mathcal M}V_Z^\top H_ZV_Z,
 \qquad \text{ with } \ H_Z=H_Z^\top,\ \ \operatorname{tr}H_Z=0,
\end{equation} 
while the quadratic form $q_r$ has at most one negative eigenvalue.
\end{proposition}

\begin{proof}
We first consider Cases~(i) and~(ii).
Let us fix $c=2\T_{\rm mix}(P)/S$ in the definition of $J$.

Locally, choosing consistent branches of the polar angle,
we regard the opening angle $\theta$ as a function of the vertices
$C$ and $D$ through the identity
\begin{equation}\label{cc:eq:angle}
 \theta(C,D)=\angle COD = \arg(D-B)-\arg(C-A)\,.
\end{equation}

Stationarity of $\Phi$ and~\eqref{cc:eq:normalization} in Lemma \ref{cc:lem:shape} give
\begin{equation}\label{cc:eq:stationarity-relations}
 \Rcal'_P[V]=\mu\,\theta'[V],\qquad
 J'_P[V]=S^2\mu\,\theta'[V],
\end{equation}
where
\begin{equation}\label{cc:eq:mu}
 \mu:=\frac{\mathcal F'(\theta)-2\theta\Rcal(P)}{\theta^2}.
\end{equation}
Here and below, we write $\theta'[V]$ and
$\theta''[V,V]$ for the first and second derivatives
with respect to $\varepsilon$, evaluated at $\varepsilon=0$, of the map $\varepsilon \mapsto 
\theta(C+\varepsilon V_C,D+\varepsilon V_D)$;
primes on $\mathcal F$ denote ordinary derivatives, whereas those
on $J$ and $\Rcal$ denote shape derivatives.

Differentiating~\eqref{cc:eq:Phi} twice and using
\eqref{cc:eq:normalization} in Lemma~\ref{cc:lem:shape} and the first equality in~\eqref{cc:eq:stationarity-relations}, we obtain, with
$x=S'_P[V]$ and $y=\theta'[V]$,
\[
 S^2\Phi''_{(P,T)}[V,V]
 =\theta^2J''_P[V,V]-\theta^2S^2\mu\,\theta''[V,V]+q_r(V),
\]
where $q_r(V)$ is the ``remainder'' quadratic form
\begin{equation}\label{cc:eq:qr}
 q_r(V)=-2\theta^2\Rcal(P)x^2-4\theta^2S\mu xy
 +S^2\bigl(4\theta\mu+2\Rcal(P)-\mathcal F''(\theta)\bigr)y^2.
\end{equation}

Recall that, according to~\eqref{cc:eq:polygon-hessian} in Lemma~\ref{cc:lem:shape}, we have 
\begin{equation}\label{cc:eq:polygon-hessian1} 
 J''_P[V,V]={}2\int_P|\nabla\psi|^2\,\dd x
 +\sum_{e\in\mathcal E_N}[fV_tV_\nu]_{\partial e}\\
 -2\int_{\partial P\setminus[A,B]}fV_t(V_\nu)_s\,\dd s.
\end{equation}

We use the boundary hat functions $\eta_Z$ on $\partial P$,
defined as in the proof of Proposition~\ref{cc:prop:bulk-corner}.
We observe that, for any edge $e$ joining two free vertices,  we have 
\begin{equation}\label{f:either} 
 \int_e f\eta_Z\,\dd s=0 \,, \qquad \text{  for each endpoint $Z$ of $e$}, 
\end{equation} 
where $f$ is the function defined in \eqref{f:defvarie}.  

 Indeed, 
put $Y=X$ in Case~(i) and $Y=D$ in Case~(ii), and denote the
outward normal to $[U,V]$ by $\nu_{[U,V]}$. From~\eqref{cc:eq:first} in Lemma~\ref{cc:lem:shape} and
the second equality in~\eqref{cc:eq:stationarity-relations},
stationarity with respect to $C$ yields
\[
 \Big (\int_{[A,C]}f\eta_C\,\dd s\Big )\nu_{[A,C]}
 +\Big (\int_{[C,Y]}f\eta_C\,\dd s\Big )\nu_{[C,Y]}
 =S^2\mu\,\nabla_C\theta.
\]
Since $\nabla_C\theta$ is normal to $[A,C]$ and
$\nu_{[A,C]},\nu_{[C,Y]}$ are linearly independent, we infer that 
$\int_{[C,Y]}f\eta_C\,\dd s=0$.
The analogous argument at $D$ gives the corresponding identity
on its inner edge. In Case~(i), stationarity at $X$ and
$\nabla_X\theta=0$ give both identities at $X$.
Thus, \eqref{f:either}  is fulfilled on each edge joining two free vertices. 

On such an edge, $V_t$ is a linear combination of the two endpoint
hat functions, while $(V_\nu)_s$ is constant. Therefore, the contribution of such an edge to 
the last integral in~\eqref{cc:eq:polygon-hessian1}
vanishes.  On the remaining Neumann
edges, only one endpoint is free. Their contributions to 
the last integral in~\eqref{cc:eq:polygon-hessian1}, as well as each term in the sum of the  endpoint terms in~\eqref{cc:eq:polygon-hessian1},  is of the kind 
$(V_Z\cdot t_e)(V_Z\cdot\nu_e)$.
Since $t_e\perp\nu_e$, each such term is a quadratic form in $V_Z$, which is represented by  a symmetric
$2\times2$ matrix of trace zero. Hence
\[
 J''_P[V,V]=2\int_P|\nabla\psi|^2\,\dd x+G(V),
 \qquad
 G(V)=\sum_{Z\in\mathcal M}V_Z^\top N_ZV_Z,
 \quad\operatorname{tr}N_Z=0.
\]
We now observe that  $\theta''[V,V]$ has the same structure. Indeed, 
the Hessian of $\theta$ with respect to $(C,X,D)$ in Case~(i), and with respect  to 
 $(C,D)$  in Case~(ii), is given by  
 \[
\begin{aligned} 
D^2_{(C,X,D)}\theta& =
 \begin{pmatrix}
 -D^2\arg(C-A)&0&0\\
 0&0&0\\
 0&0&D^2\arg(D-B)
 \end{pmatrix},
\\  \noalign{\bigskip} 
 D^2_{(C,D)}\theta & =
 \begin{pmatrix}
 -D^2\arg(C-A)&0\\
 0&D^2\arg(D-B)
 \end{pmatrix},
\end{aligned}
\] 
where each entry is a $2\times2$ block.  
Since $\arg$ is locally harmonic away from the origin, every diagonal block
has zero trace.

Therefore we have proved that  the representation formula \eqref{cc:eq:weighted-hessian}  
is satisfied, by taking $q_ r ( V)$ as in \eqref{cc:eq:qr}, and 
\[
 q_0(V):=\theta^2G(V)-\theta^2S^2\mu\,\theta''[V,V] \,,
\]
which satisfies the required condition \eqref{f:q0}. 
It remains to show that  $q_r(V)$ has at most one negative eigenvalue.

Completing the square in~\eqref{cc:eq:qr}, we can rewrite $q_ r ( V)$ as 
\begin{equation}\label{cc:eq:completed-square}
 q_r(V)=-2\theta^2\Rcal(P)
 \Big (x+\frac{S\mu}{\Rcal(P)}y\Big )^2
 +S^2\Big (\frac{2(\theta\mu+\Rcal(P))^2}{\Rcal(P)}
 -\mathcal F''(\theta)\Big)y^2.
\end{equation}

Since $\mathcal R(P)>0$ and $\mathcal F''(\theta)<0$,
the coefficient of $y^2$ is positive. Thus $q_r$ is
nonnegative on the subspace defined by
$
x+\frac{S\mu}{\mathcal R(P)}y=0,
$
which has codimension at most one in the space of
vertex velocities. Consequently, $q_r$ has at most
one negative eigenvalue.

In Case~(iii), we restrict to variations of $(X,D)$ at fixed $\theta$.
Stationarity of $\Phi$ then gives $\Rcal'_P[V]=0$, and the first
identity in~\eqref{cc:eq:normalization} of Lemma~\ref{cc:lem:shape} yields $J'_P[V]=0$.
In particular, $J$ is stationary under independent variations of
$X$ and $D$. Set $c=2\T_{\rm mix}(P)/S$.
Equation~\eqref{cc:eq:corner-hessian} in Proposition~\ref{cc:prop:bulk-corner} therefore gives
\[
 J''_P[V,V]=2\int_P|\nabla\psi|^2\,\dd x
 +\sum_{Z\in\{X,D\}}(2w_P(Z)-c)\,\mathcal Q_Z(V_Z).
\]
Since $\theta$ is fixed, differentiating~\eqref{cc:eq:Phi} twice
and using the second equality in~\eqref{cc:eq:normalization} of
Lemma~\ref{cc:lem:shape}, together with $J'_P[V]=0$, yields
\[
 S^2\Phi''_{(P,T)}[V,V]
 =\theta^2J''_P[V,V]
 -\frac{\theta^2c}{S}\bigl(S'_P[V]\bigr)^2.
\]
Substituting the expression for $J''_P[V,V]$ above, we obtain
\[
\begin{aligned}
 S^2\Phi''_{(P,T)}[V,V]
 ={}&2\theta^2\int_P|\nabla\psi|^2\,\dd x\\
 &+\theta^2\sum_{Z\in\{X,D\}}(2w_P(Z)-c)\,\mathcal Q_Z(V_Z)
 -\frac{\theta^2c}{S}\bigl(S'_P[V]\bigr)^2.
\end{aligned}
\]
Thus the representation~\eqref{cc:eq:weighted-hessian} holds with
\[
 q_0(V)=\theta^2\sum_{Z\in\{X,D\}}(2w_P(Z)-c)\,\mathcal Q_Z(V_Z),
 \qquad q_r(V)=-\frac{\theta^2c}{S}\bigl(S'_P[V]\bigr)^2.
\]
By~\eqref{cc:eq:corner-matrix} in Proposition~\ref{cc:prop:bulk-corner}, each block of $q_0$ is symmetric
and has trace zero. The form $q_r$ has rank at most one,
since $S'_P[V]$ is linear in $V$, and therefore has at most
one negative eigenvalue. This completes the proof for Case~(iii).
\end{proof}

\enlargethispage{3\baselineskip}

\begin{theorem}[Reduction to triangular caps]\label{cc:thm:no-nontriangular}
The functional $\Phi$ has no local maximizer $(P,T)$,
with respect to the free Euclidean coordinates introduced above,
such that one of the following holds:
\begin{enumerate}[label=(\roman*)]
\item $P$ is a pentagonal cap;
\item $P$ is a two-contact quadrilateral cap;
\item $P$ is a one-contact quadrilateral cap.
\end{enumerate}
\end{theorem}

\begin{proof}
Suppose that $\Phi$ has a local maximum in one of the three cases.
We use all free vertex variations in Cases~(i) and~(ii), and
variations of $(X,D)$ at fixed $\theta$ in Case~(iii).
Let $m$ be the number of free vertices: $m=3$ in Case~(i)
and $m=2$ in Cases~(ii) and~(iii). In each case, the
velocity space has dimension $2m$ and Proposition~\ref{cc:prop:stationary-hessian}
applies.

Each symmetric $2\times2$ block of $q_0$ has zero trace,
so $q_0$ has at most $m$ negative eigenvalues. By subadditivity
of the negative index, $q_0+q_r$ has at most $m+1<2m$
negative eigenvalues. Consequently, there is a nonzero
vertex velocity field $V$ such that $(q_0+q_r)(V)\geq0$.
By Lemma~\ref{cc:lem:strict-energy}, the integral
in~\eqref{cc:eq:weighted-hessian} is strictly positive.
Hence $\Phi''_{(P,T)}[V,V]>0$, contradicting local maximality.
\end{proof}

\subsection{\texorpdfstring{Optimization among triangular caps}{Optimization among triangular caps}}\label{cc:sec:benchmark}\label{cc:app:triangles}
\leavevmode\par

\begin{theorem}[Sharp comparison for triangular caps]\label{cc:thm:triangular-comparison}
 Let   $\mathcal A_\theta$ denote the class of  triangles with a fixed interior angle $\theta \in (0, \pi)$ 
 at the Neumann-Neumann vertex (i.e., the Dirichlet portion of the boundary is  the side opposite to $\theta$). 

\smallskip
\noindent (i)  The maximizers of the functional $\mathcal R$ over $\mathcal A_\theta$ are precisely the isosceles triangles.

\smallskip 
\noindent (ii)  If $T$ is a fixed triangle in $\mathcal A _\theta$, every triangular cap $K\in\mathcal C_\theta(T)$ satisfies
\begin{equation}\label{cc:eq:residual-triangle}
 \Rcal(K)\le\Rcal(I_\theta),
\end{equation}
with equality if and only if $K=T$ and $T$ is similar to $I_\theta$.
\end{theorem}

\begin{proof}
(i) For $0<\theta\le\pi/2$, the result is proved in ~\cite[Theorem 7]{BBF}.
It remains to prove it for $\pi/2<\theta<\pi$.
For $a>0$, let $g(a)$ be the mixed torsional rigidity of the right
triangle $\{0<x<1,\ 0<y<ax\}$, with Dirichlet condition on $x=1$.
By~\cite[Theorem~3.1]{GHL}, the map
$\gamma\mapsto g(\tan\gamma)-\frac18\tan\gamma$ is strictly concave
on $(0,\pi/2)$, and $g(a)\ge a/16$.

Let $h$ be the altitude of $K$ from its vertex of angle $\theta$, splitting the apex
angle into $\alpha,\beta\in(0,\pi/2)$ with $\alpha+\beta=\theta$. Set
\[
m:=\tan\frac\theta2,\qquad
 t:=\frac{\tan\alpha+\tan\beta}{2m}\ge1.
\]
Then $|K|=h^2mt$. Relaxing the continuity constraint across the altitude
in the primal formulation of the mixed torsion~\eqref{eq:mixed-primal}
allows us to treat the two right triangles independently. Each  of them is a dilation by a factor $h$
of the reference triangle defining $g$, so the corresponding rigidities scale as $h^4$. Thus
\[
 \Rcal(K)\le\frac{g(\tan\alpha)+g(\tan\beta)}{m^2t^2}.
\]
By the concavity of the map $\gamma\mapsto g(\tan\gamma)-\frac18\tan\gamma$,
\[
 g(\tan\alpha)+g(\tan\beta)
 \le 2g(m)+\frac m4(t-1)
 \le 2t^2g(m),
\]
where the last inequality follows from $g(m)\ge m/16$ and
$t^2-1\ge2(t-1)$. Thus
\[
 \Rcal(K)\le\frac{2g(m)}{m^2}=\Rcal(I_\theta),
\]
where the last equality follows by reflection across the altitude of $I _\theta$.
Equality holds if and only if $\alpha=\beta$, namely when $K$ is isosceles.

\medskip 
(ii) Write $K=\triangle(A,P,B)$ with
$P\in\overline T\setminus[A,B]$. 
Then $\gamma:=\angle APB\ge\theta$, with equality if and only if $P=O$. 
By applying statement (i) above and Corollary \ref{c:monobench}, we have
$$\mathcal R (K) \leq \mathcal R (\gamma) \leq \mathcal R (\theta)\,.$$  
Equality holds precisely when
$P=O$ and $T$ is similar to $I_\theta$.
\end{proof}

\subsection{\texorpdfstring{Conclusion of the proof of Theorem \ref{thm:convex-cap-comparison}}{Conclusion of the convex-cap comparison}}\label{cc:sec:completion}
Recall the functional $\Phi$ from~\eqref{cc:eq:Phi} in
Subsection~\ref{cc:sec:exclusion}, defined for any
$Q\in\mathcal C_\theta(T)$ by the formula:
\[
\Phi(Q,T)=\theta^2\Rcal(Q)-\mathcal F(\theta)
=\theta^2\bigl(\Rcal(Q)-\Rcal(I_\theta)\bigr).
\]
Thus Theorem~\ref{thm:convex-cap-comparison} amounts to proving
$\Phi(Q,T)\le0$, with equality precisely when $Q=T$
and $T$ is isosceles.

Suppose by contradiction that $\Phi(Q,T)>0$ for some cap $Q$ and containing
triangle $T$. Corollary~\ref{cc:cor:large-side-ratio} excludes
side ratios $t\ge2$. Since, by
Proposition~\ref{cc:prop:bounds},
$\Rcal(Q)\le\frac16\cot\frac\theta2$, we have
\begin{equation}\label{cc:eq:angular-endpoints}
\Phi(Q,T)\le\theta^2\Rcal(Q)
\le\frac{\theta^2}{6}\cot\frac\theta2\longrightarrow0
\quad\text{as }\theta\downarrow0\text{ or }\theta\uparrow\pi.
\end{equation}
Normalize the base by a rigid motion and dilation to
$A=(0,0)$, $B=(1,0)$. For a maximizing sequence with positive
limit supremum,~\eqref{cc:eq:angular-endpoints} confines $\theta$
to a compact subinterval of $(0,\pi)$; Corollary~\ref{cc:cor:large-side-ratio}
forces $1\le t<2$. These angle and side-ratio bounds keep the
normalized apices in a compact nondegenerate family.
Moreover,~\eqref{cc:eq:uniform-angle} gives
$\Phi(Q,T)\le\theta^2\delta\cot(\theta/2)/6$; hence the
relative cap areas $\delta=|Q|/|T|$ cannot vanish along a
positive maximizing sequence. By Theorem~\ref{cc:thm:five}, we may
take five-vertex caps, and their compact parametrization together
with Lemma~\ref{cc:lem:cap-continuity} yields a
positive-area maximizer $(P,T)$. The strict comparison at $t=2$
excludes this boundary value, so $t<2$. The condition $t\ge1$
is automatic and does not restrict local variations.

If $P$ is not triangular, it has one of the three configurations
of Subsection~\ref{cc:sec:exclusion}. Since the maximum is taken
over both caps and containing triangles, it is a local maximum
under the variations of Theorem~\ref{cc:thm:no-nontriangular}
(in Case~(iii), even the variations with $\theta$ fixed suffice).
This is impossible. If $P$ is triangular,
Theorem~\ref{cc:thm:triangular-comparison}(ii) gives $\Phi(P,T)\le0$,
again a contradiction. This proves inequality~\eqref{eq:convex-cap-comparison}.

Now suppose equality holds for $Q\in\mathcal C_\theta(T)$.
By Theorem~\ref{cc:thm:five}, there is
$P\in\mathcal P_\theta(T)$ with $|P|=|Q|$ and
$\Rcal(P)\ge\Rcal(Q)=\Rcal(\theta)$.
Since $\Phi\le0$ for every admissible pair, it follows
that $\Phi(P,T)=0$. Hence $(P,T)$ is a global maximizer of $\Phi$.
Theorem~\ref{cc:thm:no-nontriangular} excludes every nontriangular
configuration, so $P$ is triangular. By
Theorem~\ref{cc:thm:triangular-comparison}(ii), $P=T$ and $T$ is
isosceles. Since $Q\subseteq T$ and $|Q|=|P|=|T|$, we also
have $Q=T$.

\section{Proof of Theorem \ref{thm:sectorial-comparison}} \label{sec:spread}

Fix $j$ and consider the sectorial pair $(Q,\Delta):=(Q_j,\Delta_j)$, with opening 
$\theta$ and base $\Gamma_\theta= [A,B]$. Let $U=D\cap Q$, and
$\Gamma=\partial U\cap\partial D$.

Since the sectorial pair is simple with respect to $D$,
every side of $D$ intersecting $Q$
determines a chord of $Q$ with both endpoints on $\partial Q$.
These chords do not cross; their endpoints may be $A$
or $B$, but cannot lie in $\operatorname{relint}\Gamma_\theta$.
Indeed, $\Gamma_\theta\subset\partial K$ lies on a supporting
line of $K$, which no side of $D$ crosses, and simplicity excludes
vertices on its relative interior. The noncrossing full chords cut
$Q$ into convex subregions. Each component of $U$ is one of these
subregions.
 
We first prove that every connected component $V$ of $U$
satisfies
\begin{equation}\label{eq:sectorial-component-bound}
\T_{\rm mix}(V, \partial V\cap\partial D )
\le |V|^2\Rcal(\theta).
\end{equation}
Every subregion either adjoins the entire base
$\Gamma_\theta$, which is then a Dirichlet side, or has a
Dirichlet chord $F$ separating it from the base. To verify the
geometric containment, use affine ray coordinates in which
$\overline\Delta=\{x,y\ge0:x+y\le1\}$ and
$\overline Q=\{(x,y):0\le x\le1,\ \ell(x)\le y\le1-x\}$.
Here $\ell$ is convex, nonnegative, with $\ell(1)=0$;
therefore it is nonincreasing. For a separating chord $F=[p,q]$,
write $p=(x_1,y_1)$ and $q=(x_2,y_2)$, with
$x_1<x_2$ and $y_1>y_2$. Convexity of $\ell$ shows that the
part of $Q$ opposite to the base is contained in
$\operatorname{conv}\{p,q,(x_1,y_2)\}$: outside
$[x_1,x_2]$ its lower boundary lies above the line through
$p,q$, and within that interval $y\ge y_2$. This also
covers chords ending at $A$ or $B$. Transforming back gives a
triangle $T_F$ with base $F$, whose other sides are parallel to
the two lateral sides of $\Delta$ and whose opening is $\theta$.
If the subregion adjoins $\Gamma_\theta$, use $T_F=\Delta$
and $F=\Gamma_\theta$. Thus in either case
$V\in\mathcal C_\theta(T_F)$.
 Since
$F\subset\partial V\cap\partial D$, Neumann relaxation and
Theorem~\ref{thm:convex-cap-comparison} yield
\[
\T_{\rm mix}(V,\partial V\cap\partial D )
\le \T_{\rm mix}(V,F)
=|V|^2\Rcal(V,F)
\le |V|^2\Rcal(\theta) \,.
\]

Now let $V_1,\ldots,V_r$ be the connected components
of $U$. By additivity of mixed torsional rigidity
and~\eqref{eq:sectorial-component-bound},
$$
\T_{\rm mix}(U,\Gamma)
=\sum_{h=1}^r
 \T_{\rm mix}(V_h,\partial V_h\cap\partial D)\le\Rcal(\theta)\sum_{h=1}^r|V_h|^2
\le\Rcal(\theta)|U|^2.
$$
If equality holds, the last inequality forces $r=1$,
while equality in the convex-cap comparison implies
$U=T_F$ and $T_F$ is isosceles. By strict monotonicity
under Neumann relaxation, equality also requires that
no additional Dirichlet side of positive length be present.
Indeed, equality in the presence of an additional open segment $G\subset\Gamma\setminus F$ would force the solution with Dirichlet condition only on $F$ to vanish on $G$. Even reflection across this Neumann side and the strong maximum principle contradict its positivity there.
Thus $\Gamma=F$, up to sets of zero length.
The converse follows directly from the definition
of $\Rcal(\theta)$.
\qed

\end{document}